\documentclass[a4paper,12pt,reqno]{article}
\pdfoutput=1

\usepackage{amssymb}
\usepackage{amsmath}
\usepackage{graphicx}

\usepackage[latin1]{inputenc}
\usepackage{mathrsfs} 
\usepackage{bm}
\usepackage{upgreek}
\usepackage[bookmarksnumbered,colorlinks, linkcolor=black, citecolor=black]{hyperref}
\usepackage[top=3cm, bottom=3cm, left=3cm, right=3cm]{geometry}

\newcommand\vm[1]{\bm{\mathrm{#1}}} 

 \usepackage{color}
 \newcommand{\red}[1]{{#1}}

\title{Enhanced error estimator based on a nearly equilibrated moving least squares recovery technique for FEM and XFEM}
  
\author{J.~J.~Ródenas$^{1}$ \and O.~A.~González-Estrada$^{2}$\address \and F.~J.~Fuenmayor$^{1}$ \and F.~Chinesta$^{3}$}

\newcommand{\address}[0]{\footnote{Email: estradaoag@cardiff.ac.uk}}

\begin{document}

\maketitle


\begin{center}\small{
$^{1}$Centro de Investigaci\'on de Tecnolog\'ia de Veh\'iculos(CITV), \\Universitat Polit\`{e}cnica de Val\`{e}ncia, E-46022-Valencia, Spain.\\ 
$^{2}$ Institute of Mechanics \& Advanced Materials, Cardiff University, School of Engineering, Queen's Building, The Parade, Cardiff CF24 3AA Wales, UK.\\
$^{3}$ EADS Corporate International Chair. École Centrale de Nantes, Nantes, France. \\}
\end{center}

\begin{abstract}

In this paper a new technique aimed to obtain accurate estimates of the error in energy norm using a moving least squares (MLS) recovery-based procedure is presented. We explore the capabilities of a recovery technique based on an enhanced MLS fitting, which directly provides continuous interpolated fields, to obtain estimates of the error in energy norm as an alternative to the superconvergent patch recovery (SPR). Boundary equilibrium is enforced using a \textit{nearest point} approach that modifies the MLS functional. Lagrange multipliers are used to impose a nearly exact satisfaction of the internal equilibrium equation. The numerical results show the high accuracy of the proposed error estimator. 
\end{abstract}



\section{Introduction}

During the last few decades numerical techniques, such as the finite element method (FEM), have been used to approximate the solution of real problems. In order to assess the quality of these approximations it is necessary to evaluate the error obtained in the simulation. Current methods used to estimate the discretization error of finite element (FE) solutions are usually classified into different families: residual based, recovery based, dual analysis techniques,... \cite{Ainsworth2000, Bangerth2003, Almeida2006}. The use of recovery based estimators is widespread due to their robustness and simple implementation into existing FE codes.

Today, novel techniques such as the extended finite element method (X\-FEM) are being used to introduce a priori known information about the problem solution into the FE formulation. The XFEM \cite{Moes1999} enriches the classical FEM basis functions using a partition of unity approach in order to capture the local features of the solution in a cracked domain, i.e. the discontinuity of the displacement field along the crack faces and the singularity of the stresses in the vicinity of the crack tip.

Although the XFEM provides highly accurate solutions, and significantly improves the modelling of certain types of problems, there is an urge to develop error control techniques for these kind of methods, mostly because of their increasing importance and the fact that they use rather coarse discretizations. For example, recovery based error estimators for partition of unity methods have been developed in \cite{ Bordas2007, Duflot2008, Rodenas2007a, Rodenas2010}, using the residual approach in \cite{Strouboulis2006, Pannachet2008} and the constitutive relation error (CRE) in \cite{Panetier2010}. In \cite{panetierladeveze2009} the CRE is used for goal oriented error estimation in XFEM. 


The enforcement of the internal and the boundary equilibrium equations for stress recovery has been previously considered in the literature as a mean to improve the quality of the recovered field. For patch based formulations, \cite{Wiberg1994, Blacker1994} introduced the squares of the residuals of the equilibrium equations to the least squares functional solved in the recovery process through a penalty parameter. In \cite{Kvamsdal1998} a point-wise enforcement of the internal equilibrium in the polynomial basis, used to represent the recovered stress at the support of each node, was presented. Then, boundary equilibrium conditions were applied on a set of sampling points in the part of the patch boundary that coincides with the domain boundary. The SPR-C technique proposed in \cite{Rodenas2007} imposed equilibrium constraints to the polynomial basis via Lagrange multipliers. The internal equilibrium was exactly satisfied at each patch, and a Taylor expansion of the applied stresses was enforced along the Neumann contour. Then, a conjoint polynomial procedure \cite{Blacker1994} was used to obtain a continuous stress field. Later, in \cite{Rodenas2007a} this technique was extended to XFEM approximations. 
 
Procedures to smooth or to recover the stress field based on MLS have also been used. In \cite{Tabbara1994} a continuous stress field was obtained through local interpolation of the nodal displacement values using MLS. In \cite{Bordas2007}  this same formulation was extended to XFEM problems, considering an enriched MLS basis and a diffraction criterion, and  an error estimate for enriched approximations was proposed. In \cite{Fleming1997} a procedure to smooth the stresses for the meshless element free Galerkin (EFG) method using MLS shape functions was described. In \cite{Xiao2004} a so called Statically Admissible Stress Recovery Technique (SAR) that used MLS to fit the stress at sampling points  was presented. The SAR technique comprised basis functions which consider the internal equilibrium equations and the local tractions conditions along the Neumann boundary. In \cite{Xiao2006} an extension of the SAR technique for XFEM was presented. Following the definition of \textit{pseudo-divergence-free field} used in \cite{Huerta2004}, we can say that in references \cite{Xiao2004} and \cite{Xiao2006} the authors considered a \textit{pseudo-satisfaction} of the internal equilibrium equation. They indicated that accurate stresses were obtained with SAR, but they did not go further to evaluate any error estimate. In \cite{Duflot2004} dual techniques were used in meshless methods to obtain an equilibrated dual problem using MLS shape functions that approximate Airy stress functions. 

The objective of this paper is to present an enhanced version of the MLS recovery technique to evaluate accurate estimates of the discretization error for FEM and XFEM problems that is based in the ideas presented in \cite{Rodenas2007a, Rodenas2010, Rodenas2007,Diez2007}. The rationale behind the proposed technique is to try to enforce the recovered stress field to satisfy continuity (this property is directly provided following the MLS approach) and the equilibrium equations that are satisfied by the exact solution such that recovered stresses get closer to the exact stress field. An appropriate application of the equilibrium constraints is required to avoid discontinuities in the recovered field. \red{For that reason, a novel approach to introduce the internal and boundary equilibrium constraints is proposed.} 
\red{The procedure has been implemented in a FE code where mesh refinement is based on element splitting and the use of constrain equations (Multi Point Constraints, MPC) to force $C^0$ continuity at hanging nodes. SPR requires special treatment of these nodes because it is based on the mesh topology. The use of the proposed technique is more flexible as it is not constrained by the topology of the finite element mesh. This feature is very powerful and it allows the direct use of the technique with isogeometric analysis with \emph{h}-adaptive refinement based on T-splines and in cases where the FE mesh is missing, like with meshless methods or elements with an arbitrary number of sides.}	

Reference \cite{Diez2007} showed that upper bounds of the error in energy norm can be obtained with recovery-based error estimators if the recovered stress field is statically admissible. The recovered stresses resulting from the use of the technique proposed in this paper are continuous, satisfy the contour equilibrium equation and provide a nearly exact satisfaction of the internal equilibrium equation (more accurate than the pseudo-satisfaction of the equilibrium equation used in previous works \cite{Xiao2004, Xiao2006}). Although the upper bound property is not guaranteed, the numerical results show that the proposed technique yields sharp error estimates which nearly bound the exact error.

The paper is organised as follows: in Section ~\ref{sec:ProbStatement} we present the reference problems and their approximate solutions using FEM and XFEM. Section ~\ref{sec:ErrorEstimation} deals with the main aspects of error estimation in FE approximations and the moving least squares formulation considering equilibrium conditions. Finally, numerical results are presented in Section~\ref{sec:NumExamples}, and conclusions are drawn in Section ~\ref{sec:Conclusions}.

\section{Problem Statement and Solution}
\label{sec:ProbStatement}

\subsection{ Problem statement}
 
Let us consider the 2D linear elasticity problem. The unknown displacement field $\vm{u}$, taking values in $\Omega \subset \mathbb{R}^{2}$, is the solution of the boundary value problem given by 
\begin{align}
  -\nabla \cdot \vm{\upsigma} \left(\vm{u}\right) &= \vm{b}  	&&  {\rm in }\; \Omega 
   \label{Eq:IntEq} \\
   \vm{\upsigma} \left(\vm{u} \right)\cdot \vm{n} &= \vm{t} 	&&  {\rm on }\; \Gamma _{N}  		\label{Eq:Neumann}\\
   \vm{u}                                         &= \vm{0}	&&  {\rm on }\; \Gamma _{D} \label{Eq:Dirichlet}
\end{align}

\noindent where $\Gamma _{N}$ and $\Gamma _{D}$ denote the Neumann and Dirichlet boundaries with $\partial \Omega = \Gamma_N \cup \Gamma_D$ and $\Gamma_N \cap \Gamma_D =\varnothing$. The Dirichlet boundary condition in (\ref{Eq:Dirichlet}) is taken homogeneous for the sake of simplicity. 
 
The weak form of the problem reads: Find $\vm{u} \in V$ such that 
\begin{equation} \label{Eq:WeakForm} 
\forall \vm{v} \in V \qquad a(\vm{u},\vm{v}) = l(\vm{v}),
\end{equation}  
 
\noindent where $\vm{V}$ is the standard test space for the elasticity problem such that $V = \{\vm{v} \;|\; \vm{v} \in  H^1(\Omega) , \vm{v}|_{\Gamma_D}(\vm{x}) = \vm{0} \}$, and 
\begin{align}
a(\vm{u},\vm{v})  & := \int _{\Omega}  \vm{\sigma}^{T} (\vm{u})  \vm{\varepsilon}(\vm{v}) d \Omega = 
\int _{\Omega}  \vm{\sigma}(\vm{u})^{T}  \vm{D}^{-1} \vm{\sigma} (\vm{v}) d \Omega \\
l(\vm{v})&:=\int _{\Omega} \vm{b}^{T}   \vm{v}d \Omega + \int _{\Gamma_N} \vm{t}^{T} \vm{v}d \Gamma,
\end{align}
 
\noindent where $ \vm{\sigma}$ and $\vm{\varepsilon}$ denote the stresses and strains, and  $\vm{D} $ is the elasticity matrix of the constitutive relation $\vm{\sigma}= \vm{D} \vm{\varepsilon} $.

\subsubsection{Singular problem:}

\begin{figure}[!ht]
	\centering
	\includegraphics{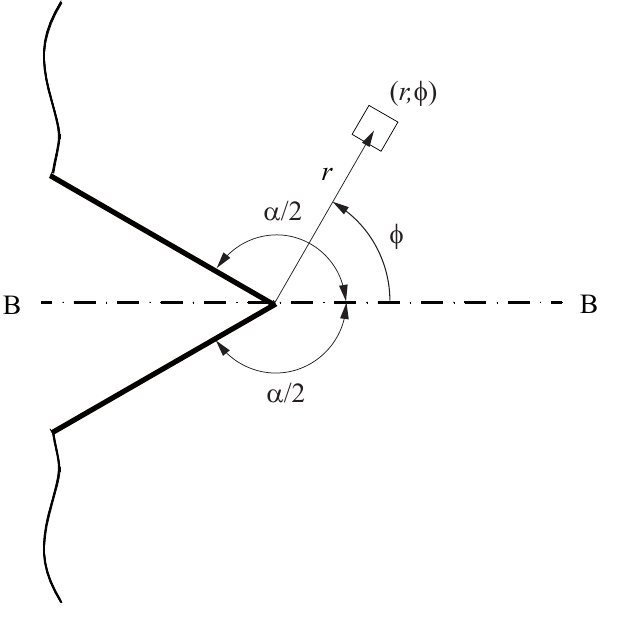} 
	\caption{Sharp reentrant corner in an infinite half-space.}
	\label{fig:Vnotch}
\end{figure}

Figure~\ref{fig:Vnotch}  shows a portion of an elastic body with a reentrant corner (or V-notch), subjected to tractions on remote boundaries. For this kind of problems, the stress field exhibits a singular behaviour at the notch vertex.

The analytical solution of the stress distribution in the vicinity of the singular point is a linear combination of singular and non-singular terms. It is often claimed that the term with a highest order of singularity dominates over the rest of terms in a sufficiently close zone surrounding the singular point. The analytical solution to this singular elastic problem in the vicinity of the singular point can be found in \cite{Williams1952, Szabo1991}. If $\alpha = 2\pi$ the problem corresponds to the classic crack problem of linear elastic fracture mechanics (LEFM) that will be used in the numerical examples. The following expressions show the first term of the asymptotic expansion of the solution for mixed mode loading conditions in a 2D cracked domain \cite{Szabo1991}:

\begin{equation}\label{Eq:CrackDisp}
\begin{aligned}
u_1(r,\phi)&=\frac{K_{\rm I}}{2\mu} \sqrt{\frac{r}{2\pi }} \cos\frac{\phi }{2} \left(\kappa-\cos \phi\right)+\frac{K_{\rm II}}{2\mu} \sqrt{\frac{r}{2\pi}} \sin\frac{\phi}{2}\left(2+\kappa+ \cos \phi \right) \\
u_2(r,\phi)&=\frac{K_{\rm I}}{2\mu} \sqrt{\frac{r}{2\pi }} \sin\frac{\phi }{2} \left(\kappa-\cos \phi\right)+\frac{K_{\rm II}}{2\mu} \sqrt{\frac{r}{2\pi}} \cos\frac{\phi}{2}\left(2-\kappa- \cos \phi \right)
\end{aligned}
\end{equation}

\begin{equation}\label{Eq:CrackStress}
\begin{aligned}
\sigma_{11}(r,\phi)&=\frac{K_{\rm I}}{\sqrt{2\pi r}}\cos\frac{\phi}{2}\left(1-\sin\frac{\phi}{2} \sin\frac{3\phi}{2}\right)-\frac{K_{\rm II}}{\sqrt{2\pi r}}\sin\frac{\phi}{2}\left( 2+\cos\frac{\phi}{2}\cos \frac{3\phi}{2} \right)\\
\sigma_{22}(r,\phi)&=\frac{K_{\rm I}}{\sqrt{2\pi r}} \cos\frac{\phi}{2} \left(1+\sin\frac{\phi}{2} \sin\frac{3\phi}{2}\right)+\frac{K_{\rm II}}{\sqrt{2\pi r}}\sin\frac{\phi}{2} \cos\frac{\phi}{2}\cos \frac{3\phi}{2}\\
\sigma_{12}(r,\phi)&=\frac{K_{\rm I}}{\sqrt{2\pi r}} \sin\frac{\phi}{2} \cos\frac{\phi}{2} \cos \frac{3\phi}{2} + \frac{K_{\rm II}}{\sqrt{2\pi r}}\cos\frac{\phi}{2} \left(1-\sin\frac{\phi}{2} \sin\frac{3\phi}{2}\right)
\end{aligned}
\end{equation}

\noindent where $K_{\rm I}$ and $K_{\rm II}$ are the generalise stress intensity factors (GSIF) for modes I and II. The GSIF are multiplicative constants that depend on the loading of the problem and linearly determine the intensity of the displacement and stress fields in the vicinity of the singular point. 

\subsection{ Solution with FEM/XFEM.}
\label{Sec:FEMSolution}

Let $\vm{u}^{h}$ be a finite element approximation to  $\vm{u}$ such that $\vm{u}^{h}(\vm{x}) = \sum_{i \in I} N_i (\vm{x}) \vm{u}_i$,  where $N_i$ represent the shape functions associated with node $i$ and $I$ is the set of all the nodes in the mesh. The solution lies in a functional space $V^{h} \subset V$ associated with a mesh of isoparametric finite elements of characteristic size $h$, and it is such that 
\begin{equation}  
\forall \vm{v} \in V^{h} \qquad a(\vm{u}^{h},\vm{v}) = l(\vm{v})
\end{equation} 

%

\begin{figure}[!ht]
	\centering
	\includegraphics{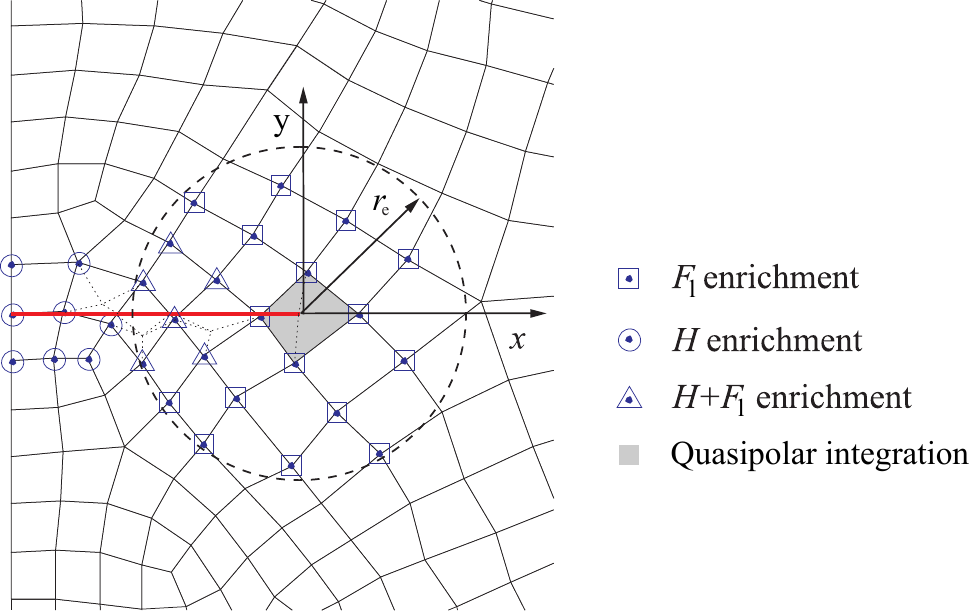} 
	\caption{Classification of nodes in XFEM. Fixed enrichment area of radius $r_e$}
	\label{fig:NodalSets}
\end{figure}

Considering an XFEM formulation for the case of the singular problems of LEFM above-mentioned, the FE approximation is enriched with Heaviside functions to describe the discontinuity of the displacement field and with crack tip functions to represent the asymptotic behaviour of the stress field near the crack tip. This avoids the need for a conforming mesh to describe the geometry of the crack \cite{Belytschko1999} and the use of adaptive techniques in order to capture the special features of the solution. To ensure the continuity of the solution, the partition of unity property of the classical linear shape functions is used. Therefore, the XFEM displacements interpolation in a 2D model is given by:
\begin{equation} \label{Eq:uXFEM} 
\vm{u}^{h}(\vm{x}) =\sum _{i\in I} N_{i}(\vm{x}) \vm{a}_{i}  +\sum _{j\in J}N_{j}(\vm{x}) H(\vm{x})\vm{b}_{j}  +\sum _{m\in M}N_{m}(\vm{x}) \left(\sum _{\ell =1}^{4}F_{\ell } (\vm{x})\vm{c}_{m}^{\ell }  \right)  
\end{equation}
 
\noindent where $\vm{a}_{i}$ are the conventional nodal degrees of freedom, $\vm{b}_{j}$ are the coefficients associated with the discontinuous enrichment functions, and $\vm{c}_{m}$ those associated with the functions spanning the asymptotic field. In the above equation, $I$ is the set of all the nodes in the mesh, $M$ is the subset of nodes enriched with crack tip functions, and $J$ is the subset of nodes enriched with the discontinuous enrichment (see Figure~\ref{fig:NodalSets}). In (\ref{Eq:uXFEM}) the Heaviside function $H$, with unitary modulus and a change of sign on the crack face, describes the displacement discontinuity if the finite element is intersected by the crack. $F_{\ell }$ are the set of branch functions used to represent the asymptotic expansion of the displacement field around the crack tip seen in (\ref{Eq:CrackDisp}). The $F_{\ell }$ functions used in this paper for the 2D case are \cite{Belytschko1999}:
\begin{equation} 
\left\{F_{\ell } \left(r,\phi \right)\right\}\equiv \sqrt{r} \left\{\sin \frac{\phi }{2} ,\cos \frac{\phi }{2} ,\sin \frac{\phi }{2} \sin \phi ,\cos \frac{\phi }{2} \sin \phi \right\}
\end{equation}

More details regarding the XFEM implementation are given in the numerical examples section.

\section{Nearly Equilibrated MLS Recovery Technique For Error Estimation In Energy Norm}
\label{sec:ErrorEstimation}

The discretization error is defined as the difference between the exact solution $\vm{u}$ and the finite element solution $\vm{u}^{h}$: $ \vm{e} = \vm{u} -\vm{u}^{h} $. Since the exact solution is in practice unknown, in general, the exact error can only be estimated. To obtain an estimation of $\vm{e}$, measures in terms of the energy norm are normally used. The Zienkiewicz-Zhu error estimator defined as 
 
\begin{equation} \label{Eq:ZZestimator} 
\lVert\vm{e}\lVert ^{2} \approx \lVert \vm{e}_{es}\lVert^{2}=\int _{\Omega}\left( \vm{\sigma}^*- \vm{\sigma}^h \right)^{T} \vm{D}^{-1} \left(\vm{\sigma}^*- \vm{\sigma}^h \right)d\Omega   
\end{equation} 
 
\noindent relies on the recovery of an improved stress field $\vm{\sigma} ^{*} $, which is supposed to be more accurate than the FE solution $\vm{\sigma} ^{h} $, to obtain an estimation of the error in energy norm $\lVert \vm{e}_{es}\lVert$.

The proposed recovery technique is based on a moving least squares procedure which provides a continuous recovered stress field \cite{liu2003} and will be denoted as MLSCX. The satisfaction of boundary and internal equilibrium has been considered in the formulation, aiming to create a statically admissible stress field that would provide accurate error estimates in the FEM and XFEM frameworks. 
 
\subsection{ MLS recovery} 
 
The MLS technique is based on a weighted least squares formulation biased towards the test point where the value of the function is asked. The technique considers a polynomial expansion for each one of the components of the recovered stress field in the form:
\begin{equation}
\sigma_{i} ^{*} (\vm{x}) = \vm{p}(\vm{x}) \vm{a}_i (\vm{x})\quad i=xx,yy,xy
\end{equation}

\noindent where $\vm{p}$ represents a polynomial basis and $\vm{a}$ are unknown coefficients
\begin{align}
\vm{p} (\vm{x}) &= \{1 \; x \; y \; x^2 \; xy \; y^2 \; \ldots \} \\
\vm{a}_i (\vm{x}) &= \{ a_{0_i}(\vm{x}) \; a_{1_i}(\vm{x}) \;a_{2_i}(\vm{x}) \;a_{3_i}(\vm{x}) \;a_{4_i}(\vm{x}) \; a_{5_i}(\vm{x}) \; \ldots \}^T
\end{align}

For 2D, the expression to evaluate the recovered stress field reads:
\begin{equation} \label{Eq:MLSGeneral}
 \vm{\sigma}^{*} (\vm{x}) = 
\begin{Bmatrix}
 \sigma_{xx}^{*}(\vm{x})\\
 \sigma_{yy}^{*}(\vm{x})\\  
 \sigma_{xy}^{*}(\vm{x})                                               
\end{Bmatrix} = 
\vm{P}(\vm{x}) \vm{A} (\vm{x}) =
\begin{bmatrix}
 \vm{p}(\vm{x}) & \vm{0} & \vm{0} \\
 \vm{0} & \vm{p}(\vm{x}) & \vm{0} \\
 \vm{0} & \vm{0} & \vm{p}(\vm{x}) 
\end{bmatrix}
\begin{Bmatrix}
 \vm{a}_{xx}(\vm{x})\\
 \vm{a}_{yy}(\vm{x})\\  
 \vm{a}_{xy}(\vm{x})  
\end{Bmatrix}
\end{equation}

The format of (\ref{Eq:MLSGeneral}), considering the three components of the stress vector in a single equation, will result useful to impose the constraints required to satisfy the equilibrium equations.

Suppose that $\vm{\chi}$ is a point within $\Omega_{\vm{x}}$, being $\Omega_{\vm{x}}$ the support corresponding to a point $\vm{x}$ defined by a distance (radius) $R_{\Omega_{\vm{x}}}$. The MLS approximation for each stress component at $\vm{\chi}$ is given by
\begin{equation}
\sigma_{i} ^{*} (\vm{x},\vm{\chi}) = \vm{p}(\vm{\chi}) \vm{a}_i (\vm{x})\quad \forall\vm{\chi}\in\Omega_{\vm{x}}, \quad i=xx,yy,xy
\end{equation}

To obtain the coefficients $\vm{A}$ we have adopted the \textit{Continuous Moving Least Squares Approximation} described in \cite{liu2003}. The following functional is minimised: 
\begin{align}  \label{Eq:MLSfunctional}
J(\vm{x})  &= \int_{\Omega_{\vm{x}}}W \left(\vm{x}-\vm{\chi} \right)\left[\vm{\sigma} ^{*} \left(\vm{x},\vm{\chi} \right)-\vm{\sigma} ^{h} \left(\vm{\chi} \right)\right]^{2}d\vm{\chi}
\end{align}

Evaluating $\partial J / \partial \vm{A} = 0$ results in the linear system $\vm{M}(\vm{x})\vm{A}(\vm{x})=\vm{G}(\vm{x})$ used to evaluate $\vm{A}$, where
\begin{equation}  \label{Eq:MLSIntegrals} 
\begin{aligned}  
\vm{M}\left(\vm{x}\right) &= \int_{\Omega_{\vm{x}}} W\left(\vm{x}-\vm{\chi} \right)\vm{P}^{T} \left(\vm{\chi} \right)\vm{P}\left(\vm{\chi} \right) d\vm{\chi} \\ \vm{G}\left(\vm{x}\right) &= \int_{\Omega_{\vm{x}}} W\left(\vm{x}-\vm{\chi} \right)\vm{P}^{T} \left(\vm{\chi} \right)\vm{\sigma} ^{h} \left(\vm{\chi} \right) d\vm{\chi}   
\end{aligned}
\end{equation}

\noindent {Assuming that there are $n$ sampling points of coordinates $\vm{\chi}_l$ $(l=1...n)$ within the support of $\vm{x}$, with weight $H_l$ and being $|\vm{J}(\vm{\chi}_{l})|$ the jacobian determinant, the expressions in (\ref{Eq:MLSfunctional}, \ref{Eq:MLSIntegrals}) can be numerically evaluated as
\begin{equation}  \label{Eq:MLSIntegrals2} 
\begin{aligned}  
J(\vm{x})  &= \sum _{l=1}^{n} W  (\vm{x}-\vm{\chi}_{l}  ) [\vm{\sigma} ^{*}  (\vm{x},\vm{\chi}_{l}  )-\vm{\sigma} ^{h}  (\vm{\chi}_{l}  ) ]^{2}  |\vm{J}(\vm{\chi}_{l})|H_l\\
\vm{M}\left(\vm{x}\right) &= 
\sum _{l=1}^{n} W\left(\vm{x}-\vm{\chi}_{l} \right)\vm{P}^{T} \left(\vm{\chi}_{l} \right)\vm{P}\left(\vm{\chi}_{l} \right)  |\vm{J}(\vm{\chi}_{l})|H_l \\ 
\vm{G}\left(\vm{x}\right) &= 
\sum _{l=1}^{n} W\left(\vm{x}-\vm{\chi}_{l} \right)\vm{P}^{T} \left(\vm{\chi}_{l} \right)\vm{\sigma} ^{h} \left(\vm{\chi}_{l} \right)  |\vm{J}(\vm{\chi}_{l})|H_l   
\end{aligned}
\end{equation}}

The integration points for the numerical evaluation of the integrals in the above equations correspond to the integration points within $\Omega_{\vm{x}}$ used in the FE analysis, for which the stress field is already available. In (\ref{Eq:MLSIntegrals})  $W$ is the MLS weighting function, which in this paper has been taken as the fourth-order spline, commonly used in the MLS related literature: 
\begin{equation} \label{Eq:WFunction}   
W (\vm{x}-\vm{\chi} ) = 
\begin{cases}
1-6s^{2} +8s^{3} -3s^{4} & {\rm if}\; |s| \le 1\\
0			 & {\rm if}\; |s| > 1
\end{cases}
\end{equation}  

\noindent where $s$ denotes the normalised distance function given by 
\begin{equation} \label{Eq:sFunction}   
s= \frac{\left\|\vm{x}-\vm{\chi}  \right\| } {R_{\Omega_{\vm{x}}} }  
\end{equation}

In the more commonly used Discrete MLS approach \cite{liu2003} the functional $J\left(\vm{x}\right)$ would be defined as:

\begin{align}  \label{Eq:MLSfunctionalDiscrete}
J\left(\vm{x}\right) = \sum _{l=1}^{n}W \left(\vm{x}-\vm{\chi}_{l} \right)\left[\vm{\sigma} ^{*} \left(\vm{x},\vm{\chi}_{l} \right)-\vm{\sigma} ^{h} \left(\vm{\chi}_{l} \right)\right]^{2}
\end{align}

This approach would thus produce similar expressions to the equations shown in (\ref{Eq:MLSIntegrals2}) with the difference that, in the continuous approach, each of the sampling points $\vm{\chi}_l$ is weighted by its associated area $ |\vm{J}(\vm{\chi}_{l})|H_l $. Our numerical experience has shown that the Continuous MLS approximation used in this paper is more accurate than the discrete approximation, especially when the distribution of the sampling points is not uniform within the support. 

Continuity in $\vm{\sigma} ^{*} $ is directly provided by the MLS procedure previously described because the weighting function $W$  ensures that stress sampling points leave or enter the support domain in a gradual and smooth manner when $x$ moves \cite{liu2003}. The following sections are devoted to the satisfaction of the equilibrium equations.

\subsection{Satisfaction of the boundary equilibrium equation} 
\label{sec:BoundaryEquil}

The boundary equilibrium equation must be satisfied at each point along the contour. In \cite{Diez2007,Rodenas2007,Rodenas2010}, where an SPR-based technique was used, the authors enforced the satisfaction of the boundary conditions in patches along the boundary using  Lagrange Multipliers to impose the appropriate constraints between the unknown coefficients to be evaluated. However, this approach produces discontinuities in a MLS formulation as we move from a support fully in the interior of the domain to a support intersecting the boundary.

\begin{figure}[!ht]
	\centering
	\includegraphics{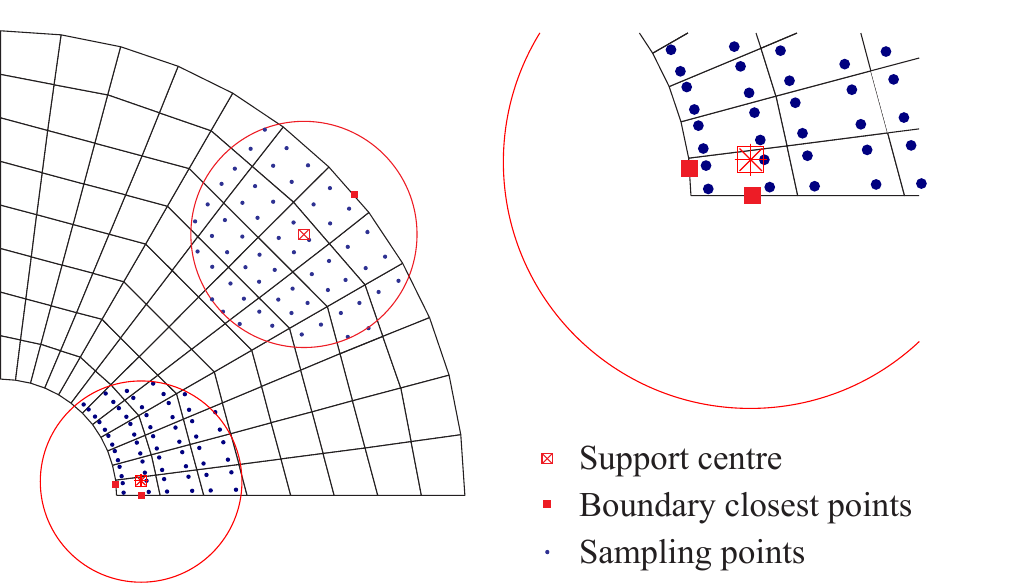}   
	\caption{MLS support with boundary conditions applied on the nearest boundary points.}
	\label{fig:MLSsupport}
\end{figure}

In order to avoid the introduction of discontinuities in the recovered field, we have followed a \textit{nearest point} approach that introduces the exact satisfaction of the boundary equilibrium equation in a smooth continuous manner. As the constraint is smoothly introduced  there is no jump when the support does not longer intersects $\Gamma$. For a point $\vm{x} \in \Omega$ whose support $\Omega_{\vm{x}}$ intersects the boundary $\Gamma$, the equilibrium constraints are considered only in the closest points $\vm{\chi}_j \in \Gamma$ on the boundaries within the support of $\vm{x}$, as shown in Figure~\ref{fig:MLSsupport}. Note that we can have more than one \textit{nearest point} for a given support, as is the case for a point $\vm{x}$ approaching a corner where we take one point for each side of the corner (see Figure ~\ref{fig:MLSsupport}). In this case, two different points have to be considered on the boundary to avoid jumps induced by the different boundary conditions when crossing the diagonal that bisects the corner.

Let us express the stress vector $\vm{\sigma}^*(\vm{x},\vm{\chi})$ in a coordinate system $\tilde{x}\tilde{y}$ aligned with the contour at $\vm{\chi}_j$ such that $\tilde{x}$ is the outward normal vector, rotated an angle $\alpha$ with respect to $x$:
\begin{equation}
 \tilde{\vm{\sigma}}^{*} (\vm{x},\vm{\chi}) = \vm{R} (\alpha) \vm{\sigma}^* (\vm{x},\vm{\chi}) 
\end{equation}  

\noindent where $\vm{R}$ is the stress rotation matrix
\begin{equation}
 \vm{R} = 
\begin{bmatrix}
 \vm{r}_{\tilde{x}\tilde{x}}\\
 \vm{r}_{\tilde{y}\tilde{y}}\\
 \vm{r}_{\tilde{x}\tilde{y}}
\end{bmatrix} =
\begin{bmatrix}
 \cos^2 \alpha 		& \sin^2 \alpha 	& \sin(2\alpha)\\
 \sin^2\alpha 		& \cos^2 \alpha 	&-\sin(2\alpha)\\
 -\sin(2\alpha)/2 	&\sin(2\alpha)/2	& \cos(2\alpha)
\end{bmatrix}
\end{equation} 

The MLS functional expressed in its continuous version and incorporating the boundary constraints reads:
\red{\begin{align}  \label{Eq:MLSfunctionalBC}
J(\vm{x})  &= 
\sum _{l=1}^{n}  W \left(\vm{x}-\vm{\chi}_{l} \right)\left[\vm{\sigma} ^{*} \left(\vm{x},\vm{\chi}_{l} \right)-\vm{\sigma} ^{h} \left(\vm{\chi}_{l} \right)\right]^{2}  |\vm{J}(\vm{\chi}_{l})|H_l  + \nonumber\\
&\quad \sum _{j=1}^{nbc} \tilde{W} \left(\vm{x}-\vm{\chi}_{j} \right) \left[\sigma^{*}_{\tilde{i}} \left(\vm{x},\vm{\chi}_{j} \right)- \sigma^{ex}_{\tilde{i}}  \left(\vm{\chi}_{j} \right)\right]^{2} \\
&= \sum _{l=1}^{n} W \left(\vm{x}-\vm{\chi}_{l} \right)\left[\vm{P}\left(\vm{\chi}_{l} \right)\vm{A}\left(\vm{x}\right)-\vm{\sigma} ^{h} \left(\vm{\chi}_{l} \right)\right]^{2}  |\vm{J}(\vm{\chi}_{l})|H_l + \nonumber\\
& \quad \sum _{j=1}^{nbc} \tilde{W} \left(\vm{x}-\vm{\chi}_{j} \right) \left[ \vm{r}_{\tilde{i}}(\alpha) \vm{P} (\vm{\chi}_{j}) \vm{A}\left(\vm{x}\right) - \sigma^{ex}_{\tilde{i}} \left(\vm{\chi}_{j} \right)\right]^{2} 
\quad \tilde{i}=\tilde{x}\tilde{x}, \tilde{x}\tilde{y} \nonumber
\end{align}
}

\noindent where $nbc$ is the number of  points $\vm{\chi}_{j}$ on the boundary where the known boundary constraints $\sigma^{ex}_{\tilde{i}}$ (in general, those would be the normal $\sigma_{\tilde{x}\tilde{x}}$ and tangential $\sigma_{\tilde{x}\tilde{y}}$ stresses) are considered. Evaluating $\partial J / \partial \vm{A} = 0$ results in the linear system $\vm{M}(\vm{x})\vm{A}(\vm{x})=\vm{G}(\vm{x})$ used to evaluate $\vm{A}$, where, in this case


\begin{align}
\vm{M} = \sum _{l=1}^{n} W\left(\vm{x}-\vm{\chi}_{l} \right)\vm{P}^{T} \left(\vm{\chi}_{l} \right)\vm{P}\left(\vm{\chi}_{l} \right)  |\vm{J}(\vm{\chi}_{l})|H_l  + \notag\\
\sum _{j=1}^{nbc} \tilde{W}(\vm{x}-\vm{\chi}_{j}) \vm{P}^{T} (\vm{\chi}_{j}) \vm{r}_{\tilde{i}}^T \vm{r}_{\tilde{i}} \vm{P} (\vm{\chi}_{j})  \\
\vm{G} = \sum _{l=1}^{n} W\left(\vm{x}-\vm{\chi}_{l} \right)\vm{P}^{T} \left(\vm{\chi}_{l} \right)\vm{\sigma} ^{h} \left(\vm{\chi}_{l} \right)  |\vm{J}(\vm{\chi}_{l})|H_l   + \notag\\
\sum _{j=1}^{nbc} \tilde{W} (\vm{x}-\vm{\chi}_{j}) \vm{P}^{T}  (\vm{\chi}_{j}) \vm{r}_{\tilde{i}}^T  \sigma^{ex}_{\tilde{i}} (\vm{\chi}_{j})
\end{align}

In the previous equations $\tilde{W}$ is a weighting function defined as:
\begin{equation} \label{Eq:WprimeFunction}
 \tilde{W} (\vm{x}-\vm{\chi}_{j} ) = 
\frac{W (\vm{x}-\vm{\chi}_{j} )}{s}=
\begin{cases}
\dfrac{1}{s} -6s +8s^{2} -3s^{3} & {\rm if}\; |s| \le 1\\
0			 	& {\rm if}\; |s| > 1
\end{cases}
\end{equation}

\begin{figure}[!ht]
  \centering
  \includegraphics[scale=0.75]{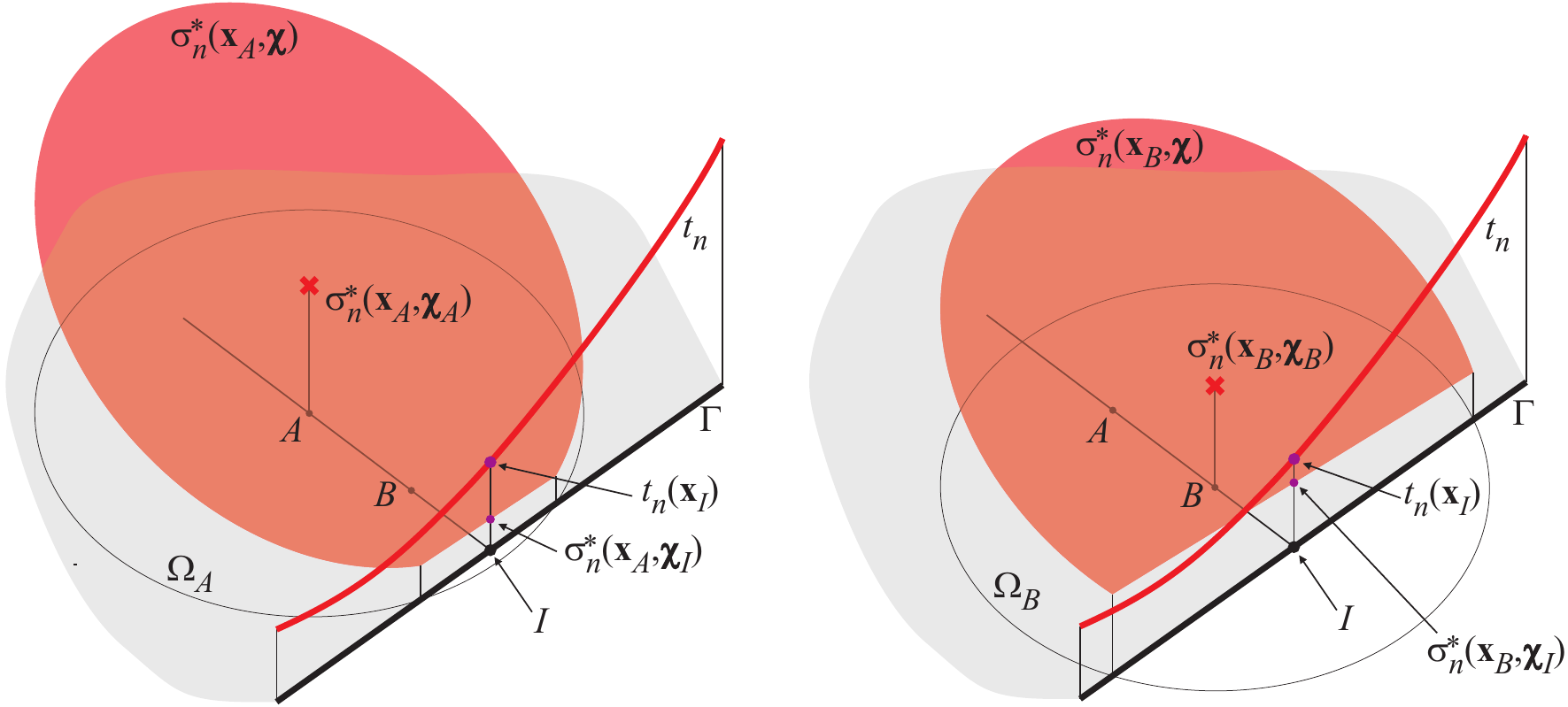}   
  \caption{Satisfaction of boundary equilibrium. $\sigma_n^*\left(\vm{x}_A,\vm{\chi}\right)$ and $\sigma_n^*(\vm{x}_B,\vm{\chi})$ are the values of $\vm{\sigma}^*\left(\vm{x},\vm{\chi}\right)$, projected along the direction normal to boundary $\Gamma$ at $I$,  in the supports $\Omega_A$ and $\Omega_B$ of the points $A$ and $B$, whose nearest point on $\Gamma$ is $I$. $t_n$ represents the normal tractions applied on $\Gamma$. Note that $\sigma_n^*\left(\vm{x},\vm{\chi}_I\right) \neq t_n\left(\vm{x}_I\right)$ although $\sigma_n^*\left(\vm{x}_B,\vm{\chi}_I\right)$ is more accurate than $\sigma_n^*\left(\vm{x}_A,\vm{\chi}_I\right)$. Thus, as $\vm{x}\rightarrow \vm{x}_I$, $\sigma_n^*\left(\vm{x},\vm{\chi}_I\right) \rightarrow t_n\left(\vm{x}_I\right) $ and, similarly the value of the stresses evaluated at the center of the support $\sigma_n^*\left(\vm{x},\vm{x}\right) \rightarrow t_n\left(\vm{x}_I\right) $ }
	\label{fig:MLSBoundary}
\end{figure}

This function has two main characteristics: 
\begin{enumerate}
  \item $\tilde{W}$ includes the weighting function $W$ such that the term for the boundary constraint is introduced smoothly into the functional $J(\vm{x})$. As a result, the recovered stress field will be continuous in $\Omega$\\

  \item $\tilde{W}$ also includes $s^{-1}$ such that the weight of the boundary constraint in $J(\vm{x})$ increases as we approach the boundary (when $\vm{x}\rightarrow \vm{\chi}_{j}$ $s \rightarrow 0$), therefore $\vm{\sigma}^*$ will tend to exactly satisfy boundary equilibrium as $\vm{x}\rightarrow \vm{\chi}_{j}$ (see Figure \ref{fig:MLSBoundary}). Note that to estimate the error using the numerical integration in (\ref{Eq:ZZestimator}), the value of $\vm{\sigma}^*$ is never evaluated on the boundary (where $s=0$) because the integration points considered are always inside the elements. 
\end{enumerate}

\subsection{Satisfaction of the internal equilibrium equation.} 

In addition to the enforcement of boundary equilibrium, we will also consider the satisfaction of the internal equilibrium equation using the Lagrange Multipliers technique. Thus, we will try to enforce the recovered stress field $\vm{\sigma} ^{*} $ to satisfy the internal equilibrium equation
\begin{equation} \label{Eq:IntEqEquat}
\nabla \cdot \vm{\sigma} ^{*} + \bm{b} = \bm{0}  
\end{equation}

 

The spatial derivatives of $\vm{\sigma} ^{*}$, considering (\ref{Eq:MLSGeneral}), are expressed as
\begin{equation}
\label{Eq:StressDeriv}\
\nabla \cdot \vm{\sigma} ^{*}=\left( \nabla \cdot \vm{P} \right) \vm{A} + \vm{P} \left( \nabla \cdot \vm{A} \right)
\end{equation}

The first terms in (\ref{Eq:StressDeriv}) can be directly evaluated differentiating the polynomial basis. Previous works \cite{Xiao2004, Xiao2006,Huerta2004} have only considered the first term in the satisfaction of the appropriate equations, thus only providing a pseudo-satisfaction of these equations \cite{Huerta2004}. Therefore, the second term in (\ref{Eq:StressDeriv}) must also be obtained. To evaluate it, we differentiate the linear system $\vm{M}\vm{A} = \vm{G}$:
%
%
\begin{equation}
\label{Eq:MAGDerivatives} 
\left( \nabla \cdot \vm{M} \right)\vm{A} + \vm{M}\left( \nabla \cdot \vm{A} \right) = 
\nabla \cdot \vm{G} 
\end{equation}

Evaluating $\nabla \cdot \vm{A}$ from (\ref{Eq:MAGDerivatives}), replacing in (\ref{Eq:StressDeriv}) and expanding leads to: 
\begin{align}    
\frac{\partial \vm{\sigma}^{*} }{\partial x} &= 
\left(\frac{\partial \vm{P}}{\partial x} -\vm{PM} ^{-1} \frac{\partial \vm{M}}{\partial x} \right)\vm{A} + \vm{PM}^{-1} \frac{\partial \vm{G}}{\partial x} = \vm{E}_{,x}\vm{A}+\vm{f}_{,x} \label{Eq:PartSigX}\\
\frac{\partial \vm{\sigma}^{*} }{\partial y} &= \left(\frac{\partial \vm{P}}{\partial y} -\vm{PM} ^{-1} \frac{\partial \vm{M}}{\partial y} \right)\vm{A} + \vm{PM}^{-1} \frac{\partial \vm{G}}{\partial y} = \vm{E}_{,y}\vm{A}+\vm{f}_{,y} \label{Eq:PartSigY}
\end{align}
     
\noindent where the partial derivatives of $\vm{M}$ and $\vm{G}$ with respect, for example, to $x$ are 
\begin{align}
\frac{\partial \vm{M}}{\partial x} &=
\sum _{l=1}^{n} \frac{\partial W (\vm{x}-\vm{\chi}_{l})}{\partial x} \vm{P}^{T}  (\vm{\chi}_{l})\vm{P} (\vm{\chi}_{l})  |\vm{J}(\vm{\chi}_{l})|H_l  + \notag\\
&\qquad \sum _{j=1}^{nbc}\frac{\partial \tilde{W} (\vm{x}-\vm{\chi}_{j})}{\partial x} \vm{P}^{T} (\vm{\chi}_{j}) \vm{r}_{\tilde{i}}^T \vm{r}_{\tilde{i}} \vm{P} (\vm{\chi}_{j})  \\
\frac{\partial \vm{G}}{\partial x} &=
\sum _{l=1}^{n} \frac{\partial W (\vm{x}-\vm{\chi}_{l})}{\partial x} \vm{P}^{T}  (\vm{\chi}_{l})\vm{\sigma} ^{h}  (\vm{\chi}_{l})   |\vm{J}(\vm{\chi}_{l})|H_l + \notag\\
&\qquad \sum _{j=1}^{nbc}\frac{\partial \tilde{W} (\vm{x}-\vm{\chi}_{j})}{\partial x} \vm{P}^{T}  (\vm{\chi}_{j}) \vm{r}_{\tilde{i}}^T     \sigma^{ex}_{\tilde{i}} (\vm{\chi}_{j})
\end{align}
 
\noindent where, differentiating (\ref{Eq:WFunction}, \ref{Eq:WprimeFunction}), 
\begin{equation}\label{Eq:dWdx}
 \frac{\partial W\left(\vm{x}-\vm{\chi} \right)}{\partial x} = 
 \frac{\partial W\left(\vm{x}-\vm{\chi} \right)}{\partial s} \frac{\partial s}{\partial x}
\end{equation}

\begin{equation}\label{Eq:dWpdx}
 \frac{\partial \tilde{W}\left(\vm{x}-\vm{\chi}_{j} \right)}{\partial x} = 
 \frac{\partial \tilde{W}\left(\vm{x}-\vm{\chi}_{j} \right)}{\partial s} \frac{\partial s}{\partial x}
\end{equation}

In these equations ${\partial s}/{\partial x}$ can be obtained from (\ref{Eq:sFunction}) or, alternatively, from (\ref{Eq:sFunctionMod}) for the case shown in the next section.  Equations (\ref{Eq:PartSigX}, \ref{Eq:PartSigY}) are expressed as a function of $\vm{A}$, so, we can write the two terms of the internal equilibrium equation (\ref{Eq:IntEqEquat}) as a function of the vector of unknowns $\vm{A}$:
\begin{align}
\frac{\partial \sigma^{*}_{xx}}{\partial x} + \frac{\partial \sigma^{*}_{xy}}{\partial y} + b_x &= \left(\vm{E}_{xx,x} + \vm{E}_{xy,y} \right) \vm{A} + \left( f_{xx,x} +  f_{xy,y} \right) + b_x = 0 \label{Eq:IntEq1} \\
\frac{\partial \sigma^{*}_{xy}}{\partial x} + \frac{\partial \sigma^{*}_{yy}}{\partial y} + b_y &=  \left(\vm{E}_{xx,y} + \vm{E}_{yy,y} \right) \vm{A} + \left(f_{xy,x} + f_{yy,y} \right) + b_y = 0 \label{Eq:IntEq2}
\end{align}

\noindent where $\Box_{i,j}$ ($i= xx, yy, xy$ and  $j=x,y$)  represents the row in $\Box_{,j}$ corresponding to the  $i^{th}$ component of the stresses. These expressions define the constraints between the coefficients $\vm{A}$ required to satisfy the internal equilibrium equation at $\vm{x}$. Lagrange Multipliers are used to impose these constraint equations.

The use of the Lagrange Multipliers technique to impose   the equilibrium constraint (\ref{Eq:IntEq1}, \ref{Eq:IntEq2}) in (\ref{Eq:MLSfunctionalBC}) leads to the following system of equations:
\begin{equation} \label{Eq:SystWithLambda}
\begin{bmatrix}
 \vm{M} & \vm{C}^T\\
 \vm{C} & \vm{0}
\end{bmatrix}
\begin{bmatrix}
 \vm{A}\\
 \vm{\lambda}
\end{bmatrix} =
\begin{bmatrix}
 \vm{G}\\
 \vm{D}
\end{bmatrix}
\end{equation} 

\noindent where $\vm{C}$ and $\vm{D}$ are the terms used to impose the constraint equations and $\vm{\lambda}$ is the vector of Lagrange Multipliers. 

However, in (\ref{Eq:MAGDerivatives}) it was assumed that $\vm{A}$ is evaluated solving $\vm{M}\vm{A} = \vm{G}$, although, operating by blocks in (\ref{Eq:SystWithLambda}) the following system of equations is obtained:

\begin{equation} \label{Eq:MACLambdaG}
\vm{M}\vm{A}+\vm{C}^T\vm{\lambda}=\vm{G}
\end{equation} 

Hence, in the formulation proposed in this paper we have neglected the term $\vm{C}^T\vm{\lambda}$  when evaluating the partial derivatives of $\vm{A}$. Evidently, this implies that the internal equilibrium equation is not fully satisfied, leading to a nearly exact satisfaction of the internal equilibrium equation. As described in the numerical examples, this approximation represents an enhancement with respect to the pseudo satisfaction of equilibrium \cite{Huerta2004}. 

References \cite{Diez2007, ladevezerougeot1999} show that the error estimator in (\ref{Eq:ZZestimator}) would produce an upper error bound if $\vm{\sigma}^*$ is statically admissible. The MLSCX recovery technique produces a continuous stress field where the internal equilibrium equation is not fully satisfied. Hence $\vm{\sigma}^*$ is continuous and nearly equilibrated and, thus, nearly statically admissible. Therefore, although the error estimate provided by the proposed recovery technique is very sharp, it is not a guaranteed upper error bound.

\subsection{ Visibility} 
 
For problems with re-entrant corners a visibility criterion is used to modify the normalised distance $s$ in (\ref{Eq:sFunction}). The standard weight function depends on the distance between the central point of the support and the sampling points, decreasing as the sampling points are located farther from the centre \cite{Bordas2007}. 

\begin{figure}[!ht]
	\centering
	\includegraphics[scale=1.2]{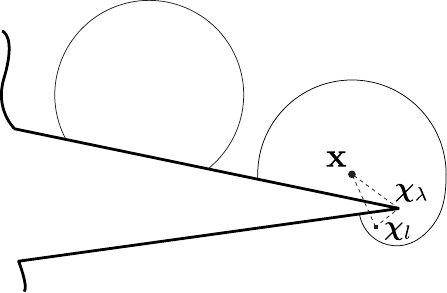}  
	\caption{Domain with re-entrant corner.}
	\label{fig:Visibility}
\end{figure}

Consider a domain with a re-entrant corner as shown in Figure~\ref{fig:Visibility}. The value of the weight function for a sampling point $\vm{\chi}_{l}$, considering a centre point $\vm{x}$ whose support contains the singularity at $\vm{\chi}_{\lambda}$, diminishes with the visibility of $\vm{\chi}_{l}$ from $\vm{x}$ such that, for points that cannot be directly viewed from $\vm{x}$, instead of (\ref{Eq:sFunction}), the following equation is used
\begin{equation}\label{Eq:sFunctionMod}
s=\frac{\left\| \vm{x}-\vm{\chi} _{\lambda } \right\| + \left\| \vm{\chi}_{l } -\vm{\chi}_{\lambda} \right\| }{R_{\Omega_{\vm{x}} }}
\end{equation}  

\subsection{ Stress splitting for singular problems.} 
\label{Sec:Splitting}
 
It is well known that smoothing techniques perform badly when the solution contains a singularity. In \cite{Rodenas2007a, Rodenas2006} a technique that decomposes the stress field in singular and smooth parts in order to improve the accuracy of SPR-based error estimators was proposed. The authors indicated that the exact stress field $\vm{\sigma}$  corresponding to a singular problem can be expressed as the contribution of a smooth stress field, $\vm{\sigma}_{smo}$, and a singular stress field, $\vm{\sigma}_{sing}$ 
\begin{equation}  
\vm{\sigma} =\vm{\sigma} _{smo} +\vm{\sigma} _{sing}   
\end{equation} 

Hence, the recovered stress field for this kind of problems can be expressed as the contribution of a smooth and a singular recovered stress fields 
\begin{equation}  
\vm{\sigma}^{*} =\vm{\sigma} _{smo}^{*} + \vm{\sigma} _{sing}^{*}   
\end{equation} 

To obtain an accurate approximation of the singular part we use the interaction integral, as shown in \cite{Shih1988}, to compute a good estimation of the GSIFs $K_{\rm I}$ and $K_{\rm II}$. Then, using the estimated values $K_{\rm I}^*$ and $K_{\rm II}^*$ we can evaluate a singular recovered stress field $\vm{\sigma} _{sing}^{*} $ from (\ref{Eq:CrackStress}). 

Assuming that $\vm{\sigma} _{sing}^{*} $ is a good approximation of the singular part $\vm{\sigma} _{sing}$, a FE-type representation of the smooth part $\vm{\sigma} _{smo}^{h} $ is given by 
\begin{equation}   
\vm{\sigma} _{smo}^{h} = \vm{\sigma} ^{h} -\vm{\sigma} _{sing}^{*}   
\end{equation} 

In \cite{Rodenas2007a, Rodenas2006} an SPR-based recovery technique was used to smooth the discontinuous stress field $ \vm{\sigma} _{smo}^{h} $. In this paper, we use the moving least squares procedure previously described to recover the smooth part of the solution $\vm{\sigma} _{smo}^{*}$. In \cite{Rodenas2007a, Rodenas2006} the stress splitting procedure was only used in a small area around the crack tip. In the procedure proposed herein the stress splitting is used in the whole domain of the problem in order to avoid discontinuities along the blending zone.  Thus, the boundary tractions to be considered for the satisfaction of the boundary equilibrium equation in the smooth problem are:
\begin{equation}
\vm{t}_{smo}=\vm{t}-\vm{t}_{sing}^*
\end{equation}

\noindent where $\vm{t}_{sing}^*$ are the projection of $\vm{\sigma} _{sing}^{*}$. It must be taken into account that the crack faces are treated as any other Neumann boundary where satisfaction of the boundary equilibrium equation will be imposed.

Note that $\vm{\sigma} _{sing}^{*}$ is equilibrated and continuous, therefore, the resulting recovered stress field $\vm{\sigma} ^{*} = \vm{\sigma} _{sing}^{*} + \vm{\sigma} _{smo}^{*}$ only has small lacks of internal equilibrium in $\vm{\sigma} _{smo}^{*}$ induced by the recovery process.

\subsection{Adaptive strategy}

\red{The refinement of the mesh using the error estimate as the guiding parameter considers an stopping criterion that checks the value of the estimated error against a prescribed or desired error. If the estimated error is higher than the desired error then the mesh is refined. Several procedures to perform the refinement are available in the literature. To define the size of the elements in the new mesh we follow the adaptive process described in \cite{ladevezemarin1992, coorevitsladeveze1995, ladevezeleguillon1983} which minimises the number of elements in the new mesh. This criterion is equivalent to the traditional approach of equally distributing the error in each element of the new mesh as proven in \cite{libettess1995, fuenmayoroliver1996}. }

%
%

\section{Numerical Examples}
\label{sec:NumExamples}

In this section numerical tests using 2D benchmark problems with exact solution are used to investigate the quality of the proposed error estimation technique. The first three problems (smooth and singular) consider a FEM approximation whilst the fourth problem is solved using an XFEM formulation. For all the models we assume a plane strain condition. Sequences of meshes with linear (TRI3), quadratic (TRI6) triangles and linear (QUAD4), quadratic (QUAD8) quadrilaterals elements are considered for the analyses. Uniform and \emph{h}-adaptive refinements have been used. The \emph{h}-adaptive refinement is based on element splitting using multipoint constraints (MPC) to impose $C^0$ continuity at hanging nodes. Quadrature rules of 1, 3, $2 \times 2$ and $3 \times 3$ Gauss points are used for TRI3, TRI6, QUAD4 and QUAD8 elements, respectively. A support size with a radius two times the average size of the surrounding elements is used to perform the MLS recovery. 19 sampling points in triangular elements and 25 sampling points in quadrilaterals are used for an accurate numerical evaluation of (\ref{Eq:ZZestimator}) in order to avoid the effect of numerical errors due to integration. The computational cost of the proposed technique could be alleviated by evaluating (\ref{Eq:ZZestimator}) using quadrature rules with fewer integration points  at the expense of introducing errors due to integration in the procedure. The MLS basis functions used in the recovery are polynomials $\vm{p}$ one order higher than the corresponding FE displacement basis. 

The performance of the technique is evaluated using the effectivity index of the error in energy norm, both at global and local levels. Globally, we consider the value of the effectivity index $\theta$ given by 
\begin{equation} \label{Eq:Effectivity}  
\theta  =\frac{\left\| \vm{e}_{es}^{} \right\| }{\left\| \vm{e}_{}^{} \right\| }   
\end{equation}  

\noindent where $\|\vm{e}_{}\| $ denotes the exact error in energy norm, and $\left\| \vm{e}_{es}^{} \right\| $ represents the evaluated error estimate. At element level, the distribution of the local effectivity index $D$, its mean value $m(|D|)$ and standard deviation $\sigma(D)$ is analysed, as described in \cite{Rodenas2007a}: 

\begin{equation} \label{Eq:LocalEffectivity}  
\begin{array}{ccc} 
{D=\theta ^{e} -1} & {\rm if} & {\theta ^{e} \ge 1} \\ 
{D=1-\dfrac{1}{\theta ^{e} }} & {\rm if} & {\theta ^{e} < 1} 
\end{array}
\qquad \qquad {\rm with} \qquad 
\theta ^{e} =\dfrac{\left\| \vm{e}_{es}^{e} \right\| }{\left\| \vm{e}_{}^{e} \right\| }   
\end{equation}  

\noindent where superscript $^e$ denotes evaluation at element level.

The h-adaptive refinement procedure considering the error in quantities of interest is implemented based on previous adaptive procedures using the error in energy norm. The technique aims to minimise the number of elements to get the target error by equally distributing the element error in the mesh.

\subsection{2$\times$2 square with a 3rd-order polynomial solution}

\begin{figure}[!ht]
  \centering
  \includegraphics{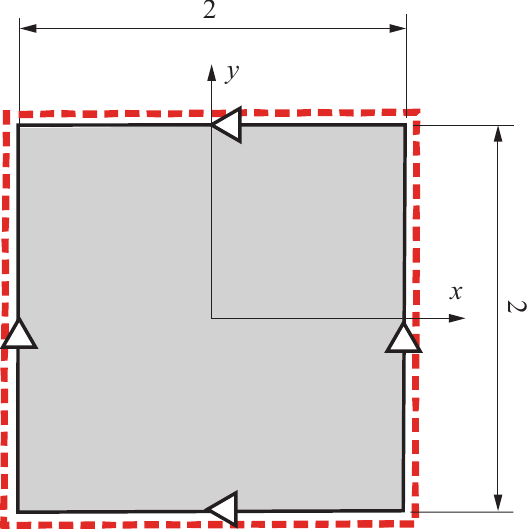}
  \caption{2$\times$2 square plate.}
  \label{fig:Square}
\end{figure}

The 2$\times$2 square model shown in Figure~\ref{fig:Square} is analysed, with material parameters $E~=~1000$ for the Young's modulus and $\nu~=~0.3$ for the Poisson's ratio. Dirichlet boundary conditions are indicated in the figure. The problem is defined such that the exact displacement solution is given by
\begin{align}
 u(x,y) &= x+x^2-2xy+x^3-3xy^2+x^2y\\
 v(x,y) &= -y-2xy+y^2-3x^2y+y^3-xy^2
\end{align}

The exact values of the stress components are applied along the Neumann boundary denoted by a dashed line in Figure~\ref{fig:Square}. These stresses can be derived from the exact displacement field under plane strain condition, and read
\begin{align}
 \sigma_{xx} &= \frac{E}{1+\nu}(1+2x-2y+3x^2-3y^2+2xy)\\
 \sigma_{yy} &= \frac{E}{1+\nu}(-1-2x+2y-3x^2+3y^2-2xy)\\
 \sigma_{xy} &= \frac{E}{1+\nu}(- x - y + \frac{x^2}{2} - \frac{y^2}{2} - 6xy)
\end{align}

The following body forces must be applied to satisfy equilibrium:
\begin{align}
 b_x(x,y) &= -\frac{E}{1+\nu}(1+y)\\
 b_y(x,y) &= -\frac{E}{1+\nu}(1-x)
\end{align}

We have used this problem to analyse the influence of different implementations of the MLS recovery technique in the error estimate, considering the following cases:

\begin{itemize}
	\item \textbf{MLS:} Plain Moving Least Squares recovery\\
	\item \textbf{MLS+BE:} MLS technique with the boundary equilibrium enhancement described in Section \ref{sec:BoundaryEquil}\\
	\item \textbf{MLS+BE+PIE:} MLS technique with the boundary equilibrium enhancement described in Section \ref{sec:BoundaryEquil} and the pseudo satisfaction of the internal equilibrium equation\\
	\item \textbf{MLSCX:} Technique proposed in this paper
\end{itemize}

\begin{figure}[ht]
      \centering
   \includegraphics{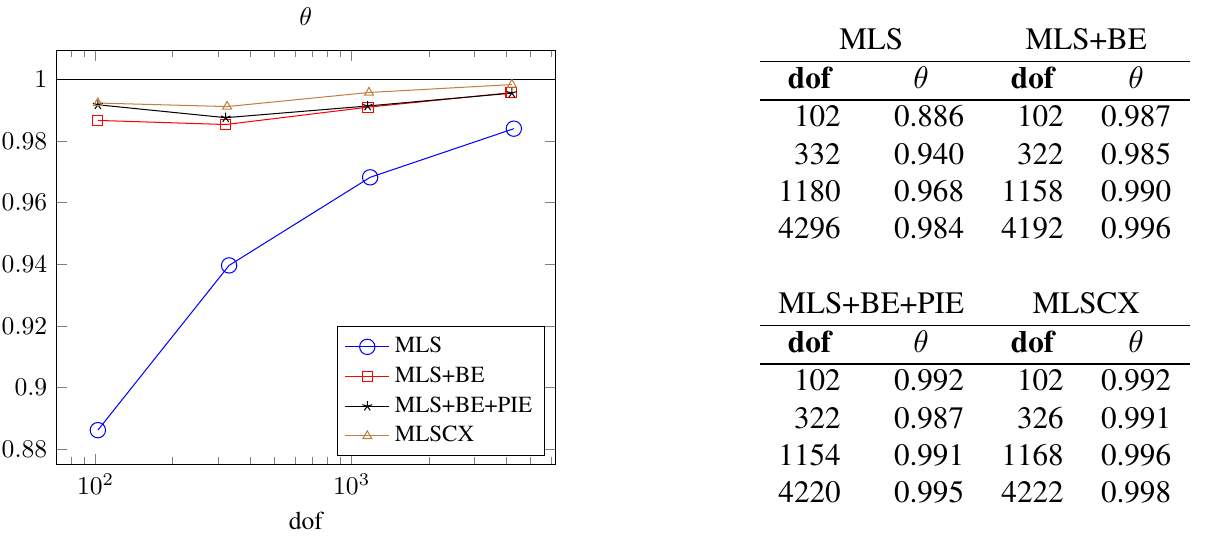}
    \caption{$2\times2$ square with \textit{h}-adaptive meshes and TRI3 elements. Evolution of $\theta$ for different recoveries}
  \label{fig:TRI3-EquilibriumSch}
\end{figure}

\begin{figure}[ht]
      \centering
   \includegraphics{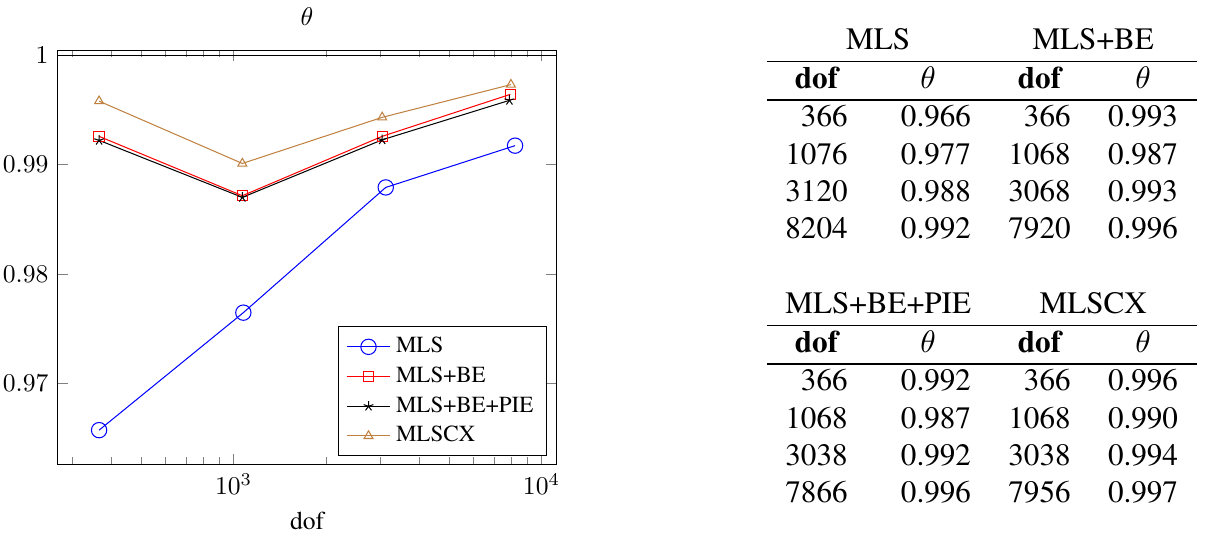}
    \caption{$2\times2$ square with \textit{h}-adaptive meshes and TRI6 elements. Evolution of $\theta$ for different recoveries.}
  \label{fig:TRI6-EquilibriumSch}
\end{figure}

\begin{figure}[ht]
      \centering
   \includegraphics{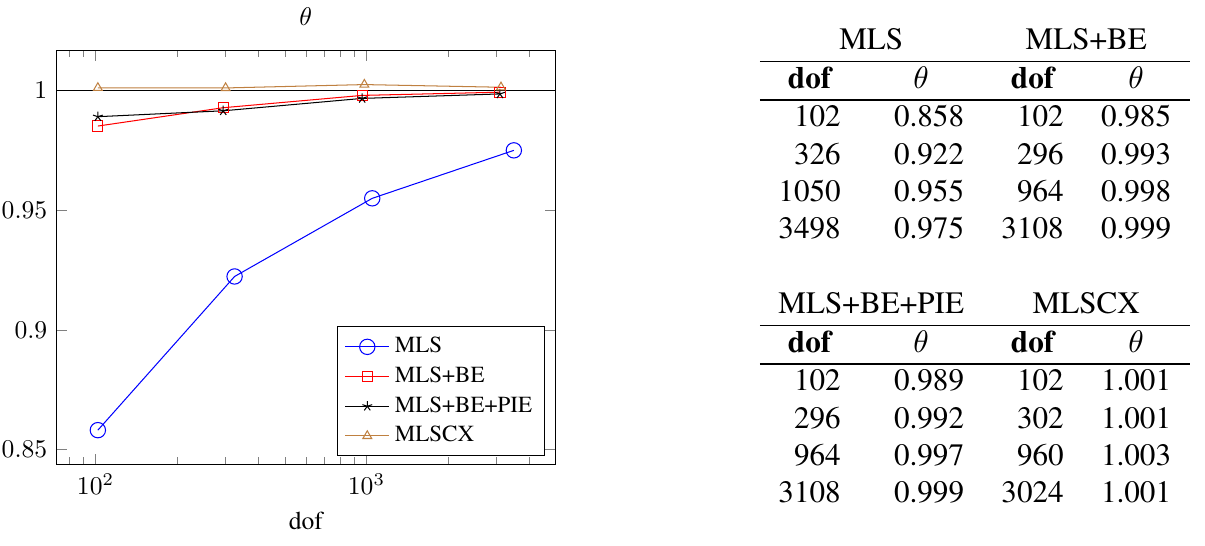}
    \caption{$2\times2$ square with \textit{h}-adaptive meshes and QUAD4 elements. Evolution of $\theta$ for different recoveries}
  \label{fig:QUAD4-EquilibriumSch}
\end{figure}

\begin{figure}[ht]
      \centering
   \includegraphics{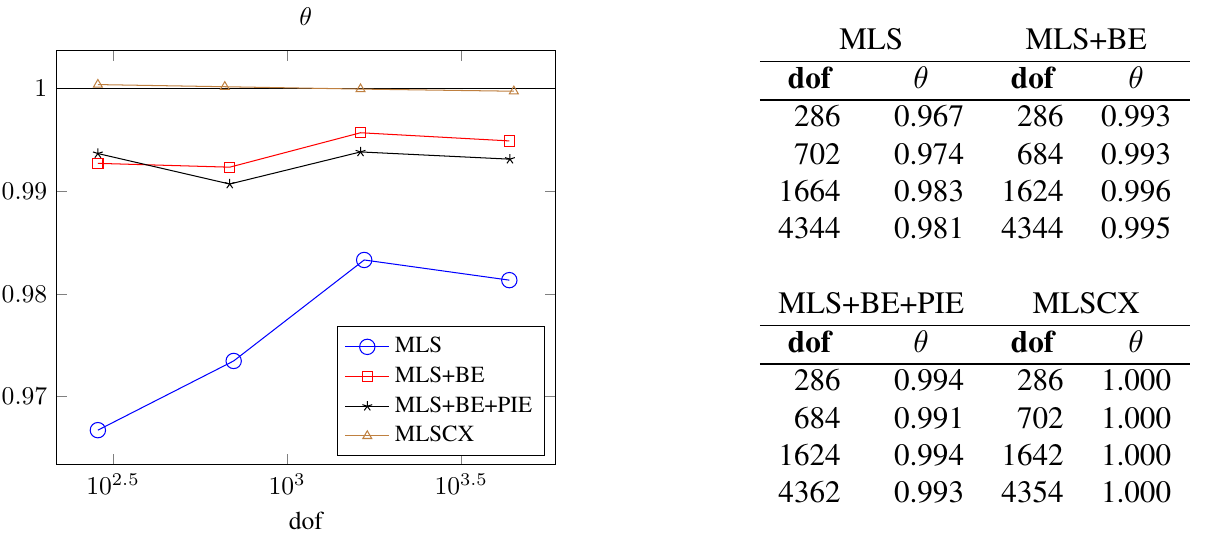}
    \caption{$2\times2$ square with \textit{h}-adaptive meshes and QUAD8 elements. Evolution of $\theta$ for different recoveries.}
  \label{fig:QUAD8-EquilibriumSch}
\end{figure}

The results for the plain MLS case will be used as reference. The other three cases represent implementations which increasingly approach the full satisfaction of the equilibrium equations. Figures \ref{fig:TRI3-EquilibriumSch} to \ref{fig:QUAD8-EquilibriumSch} show the effectivity of the error estimation vs. the number of degrees of freedom (dof) using these four implementations for \emph{h}-adaptive meshes.

These figures clearly show that the satisfaction of boundary equilibrium (curves MLS+BE) plays the most important role towards the enforcement of equilibrium and, therefore, an improvement on the accuracy of the error estimator when compared with the MLS curve. The additional pseudo-satisfaction of internal equilibrium (curves MLS+BE+PEI) does not improve, and sometimes provides worse effectivities than boundary equilibrium constraints, as it can be seen in Figure~\ref{fig:QUAD8-EquilibriumSch}. From Figures \ref{fig:TRI3-EquilibriumSch} to \ref{fig:QUAD8-EquilibriumSch} we can see an increase in the accuracy for the MLSCX curves with respect to the other curves, with effectivities very close to $\theta=1$.

\begin{figure}[!ht]
  \centering
  \includegraphics{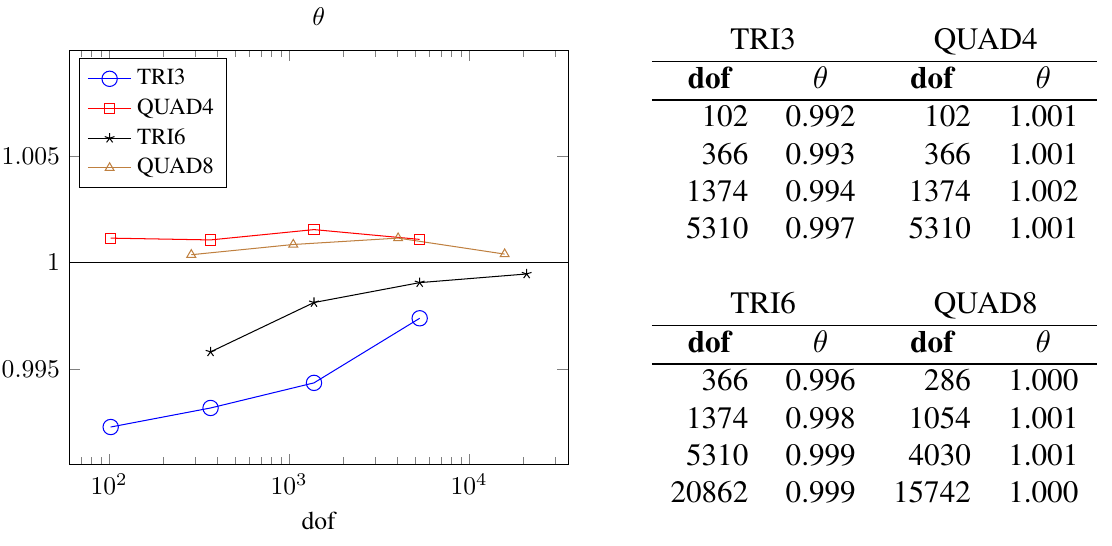}
  \caption{$2 \times 2$ square with uniformly refined meshes. Evolution of the effectivity index $\theta$ for different element types.}
  \label{fig:SquareEffectivity}
\end{figure}

Figure~\ref{fig:SquareEffectivity} shows the evolution with respect to mesh refinement of the global effectivity index $\theta$, the mean absolute value $m(|D|)$ and standard deviation $\sigma(D)$ of the local effectivity index for the different types of elements considered. Note that with the proposed technique we obtain very accurate values of $\theta$ and the error estimate converges to the exact value with the increase of the number of degrees of freedom. For this example, the best results are obtained with quadratic elements. 

\begin{figure}[!ht]
	\centering
	\includegraphics{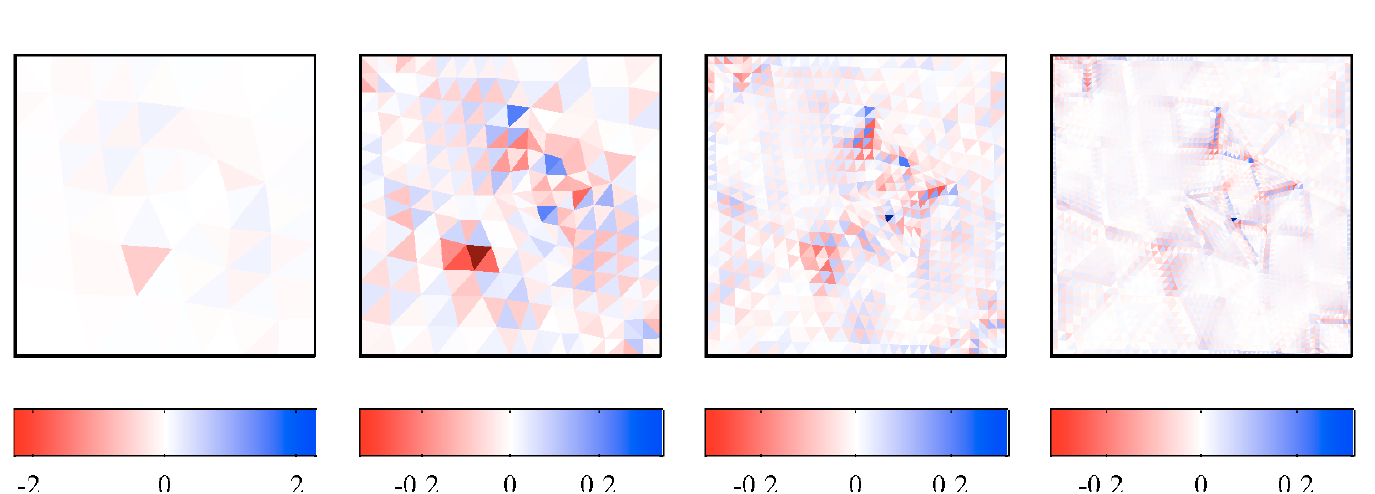}
	\caption{$2 \times 2$ square with TRI3 elements. Distribution of the effectivity index $D$ in uniformly refined meshes.}
	\label{fig:SquareDTRI3}
\end{figure}

\begin{figure}[!ht]
	\centering
	\includegraphics{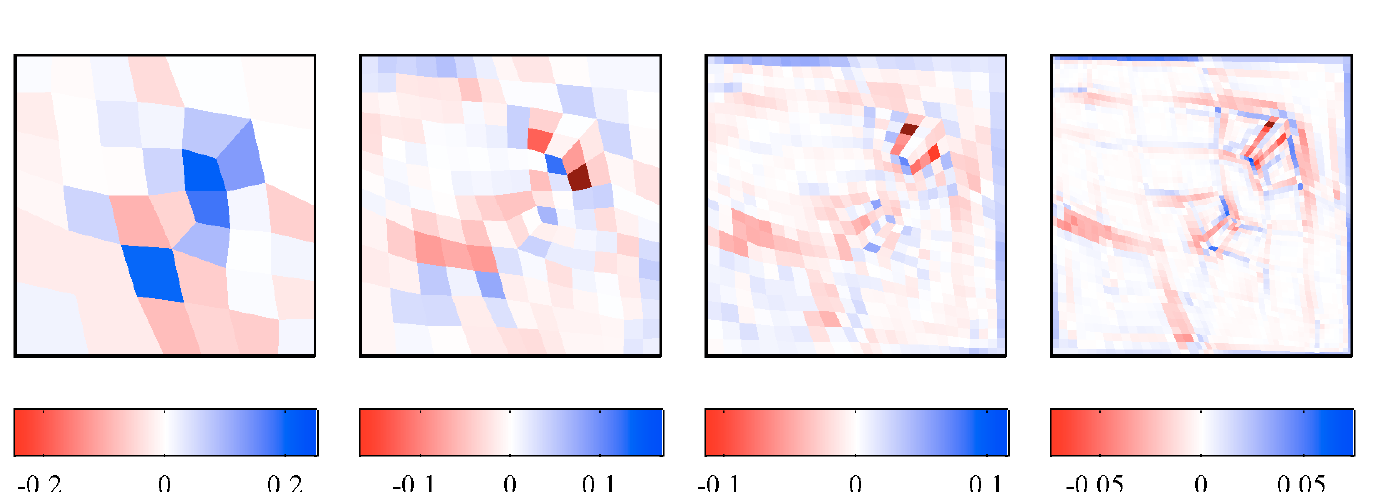}	
	\caption{$2 \times 2$ square with QUAD4 elements. Distribution of the effectivity index $D$ in uniformly refined meshes.}
	\label{fig:SquareDQUAD4}
\end{figure}

\begin{figure}[!ht]
	\centering
	\includegraphics[width=0.24\textwidth]{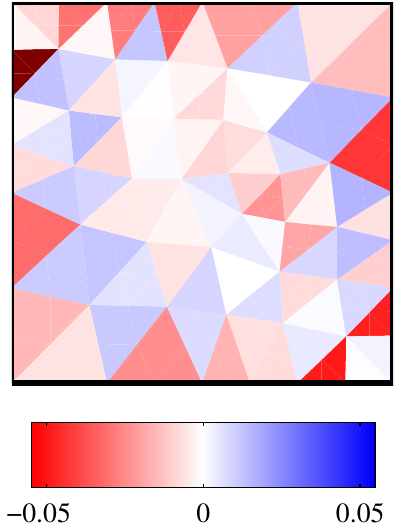}	
	\includegraphics[width=0.24\textwidth]{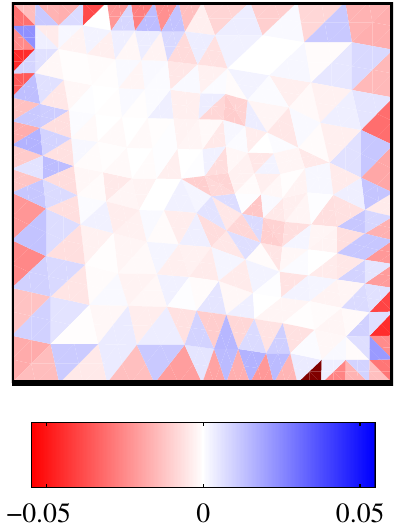}	
	\includegraphics[width=0.24\textwidth]{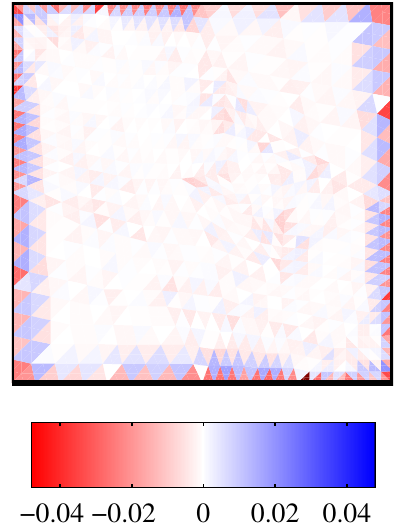}	
	\includegraphics[width=0.24\textwidth]{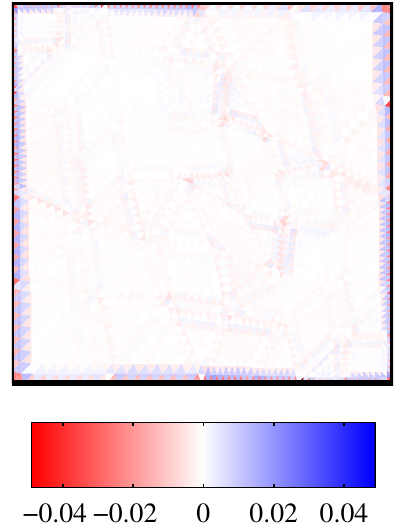}	
	\caption{\red{$2 \times 2$ square with TRI6 elements. Distribution of the effectivity index $D$ in uniformly refined meshes.}}
	\label{fig:SquareDTRI6}
\end{figure}

\begin{figure}[!ht]
	\centering
	\includegraphics[width=0.24\textwidth]{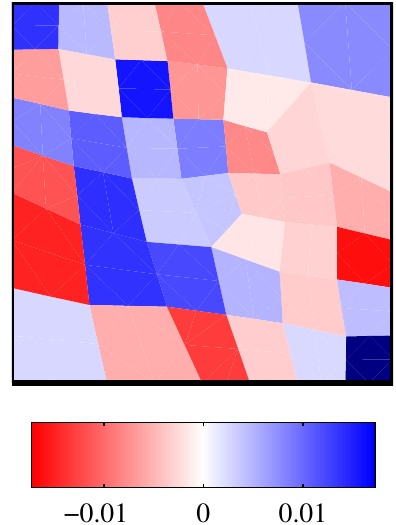}	
	\includegraphics[width=0.24\textwidth]{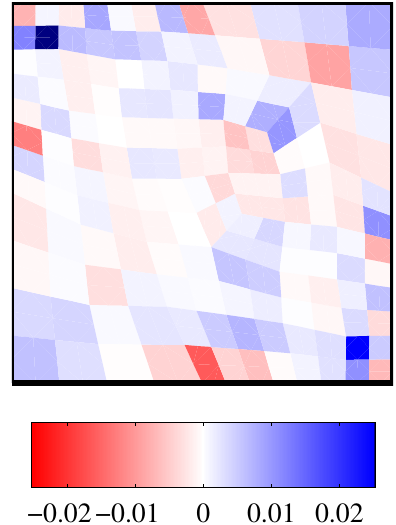}	
	\includegraphics[width=0.24\textwidth]{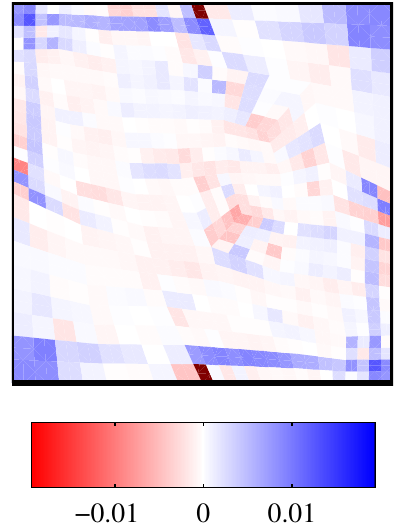}	
	\includegraphics[width=0.24\textwidth]{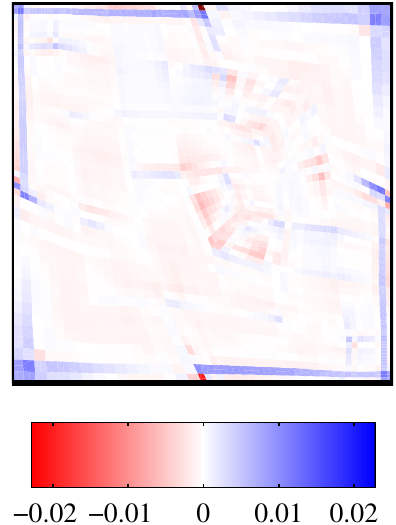}	
	\caption{\red{$2 \times 2$ square with QUAD8 elements. Distribution of the effectivity index $D$ in uniformly refined meshes.}}
	\label{fig:SquareDQUAD8}
\end{figure}

Figure~\ref{fig:SquareDTRI3} shows the distribution of the local effectivity on a set of TRI3 meshes, Figure~\ref{fig:SquareDQUAD4} displays the same results for a set of QUAD4 meshes. In both cases we can observe a quite homogeneous distribution of the local effectivity inside the domain and good results along the boundary of the problem. In addition, $D$ decreases for finer meshes and we always have  values within a very narrow range. Figures \ref{fig:SquareDTRI6} and \ref{fig:SquareDQUAD8} show a similar behaviour for quadratic elements.

\subsection{ Thick-wall cylinder subjected to an internal pressure.} 
 
\begin{figure}[!ht]
	\centering
	\includegraphics{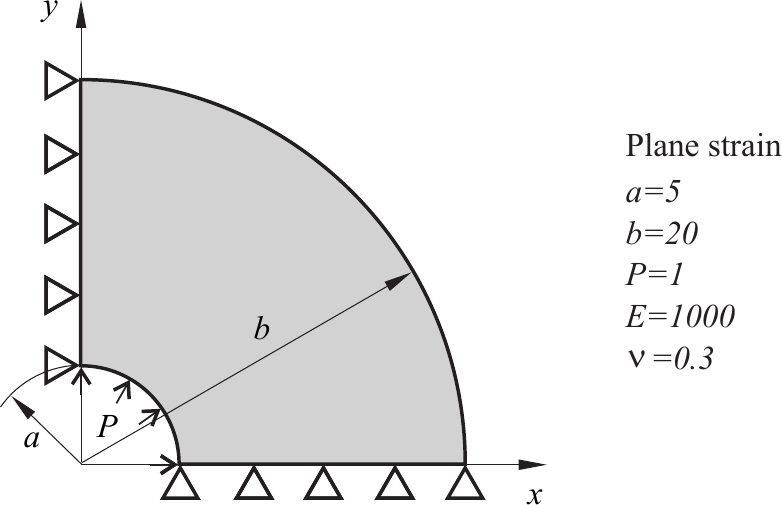}
	\caption{Thick-wall cylinder subjected to internal pressure.}
	\label{fig:cylinder}
\end{figure}

The geometrical model for this problem is shown in Figure~\ref{fig:cylinder}. Due to symmetry conditions, only one part of the section is modelled. 

The exact solution for this problem is given by the following expressions. For a point $(x,y)$, $c=b/a$, $r=\sqrt{x^{2} +y^{2} }$ the radial displacement is given by
\begin{equation} \label{Eq:uCylinder}  
u_{r} =\frac{P(1+\nu )}{E(c^{2} -1)} \left(r\left(1-2\nu \right) + \frac{b^{2} }{r} \right)  
\end{equation} 
 
Stresses in cylindrical coordinates are

 
\begin{equation} \label{Eq:stressCylinder}  
\begin{array}{cc} 
{\begin{array}{c} 
{\sigma _{r} =\dfrac{P}{c^{2} -1} \left(1-\dfrac{b^{2} }{r^{2} } \right)} \\ 
{\sigma _{t} =\dfrac{P}{c^{2} -1} \left(1+\dfrac{b^{2} }{r^{2} } \right)} 
\end{array}} & 
{\begin{array}{l} 
{\sigma _{z}  =2 \nu \dfrac{P}{c^2-1}  }
\end{array}} \end{array}  
\end{equation} 

\begin{figure}[ht]
  \centering
  \includegraphics{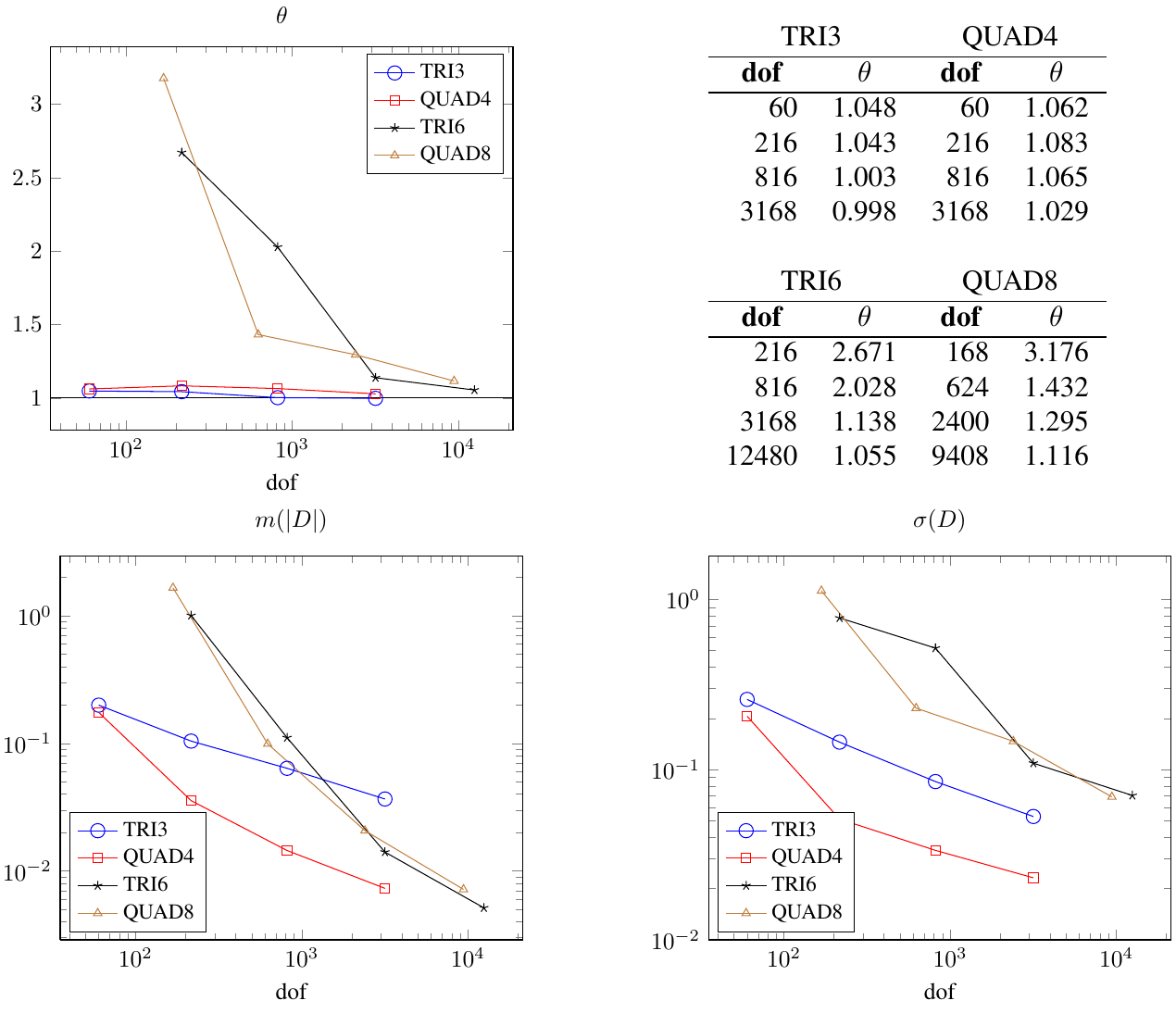}    
  \caption{Thick-wall cylinder with uniformly refined meshes. Global indicators $\theta$, $m(|D|)$ and $\sigma(D)$.}
  \label{fig:CylEffectivity}
\end{figure}

Figure \ref{fig:CylEffectivity} shows the effectivity values obtained with the MLSCX in the thick-wall cylinder using uniformly refined meshes. The results obtained are similar to those previously shown for the square plate. Compensations between underestimated and overestimated areas of the domain might result in misleading values of the global value $\theta$. To take this into account we consider the parameters $m(|D|)$ and $\sigma(D)$ which are expected to decrease when we increase the level of refinement. In Figure~\ref{fig:CylEffectivity} $m(|D|)$ and $\sigma(D)$ decrease when we increase the refinement showing a good performance of the error estimator. The results for this problem show that the proposed technique provides an accurate estimate of the exact error in energy norm. \red{On the other hand, for this example the best values are obtained when considering linear elements, thus, there is not a strong correlation between the order of the approximation and the performance of the error estimator.}

\begin{figure}[!ht]
	\centering
	\includegraphics{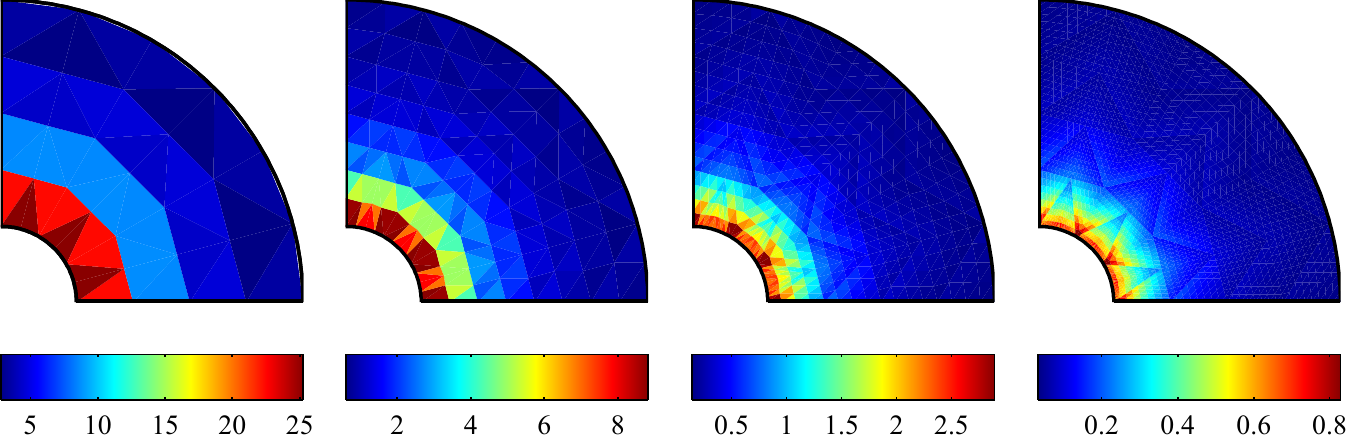}
	\caption{\red{Thick-wall cylinder with TRI3 elements. Distribution of the exact   error in energy norm in uniformly refined meshes.}}
	\label{fig:CylTRI3EnergyNorms}
\end{figure}

Figure ~\ref{fig:CylTRI3EnergyNorms} shows the distribution of the exact error in energy norm for the same meshes. The higher errors are located close to the inner radius of the cylinder. The error decreases as we increase the number of degrees of freedom.

\begin{figure}[!ht]
	\centering
	\includegraphics{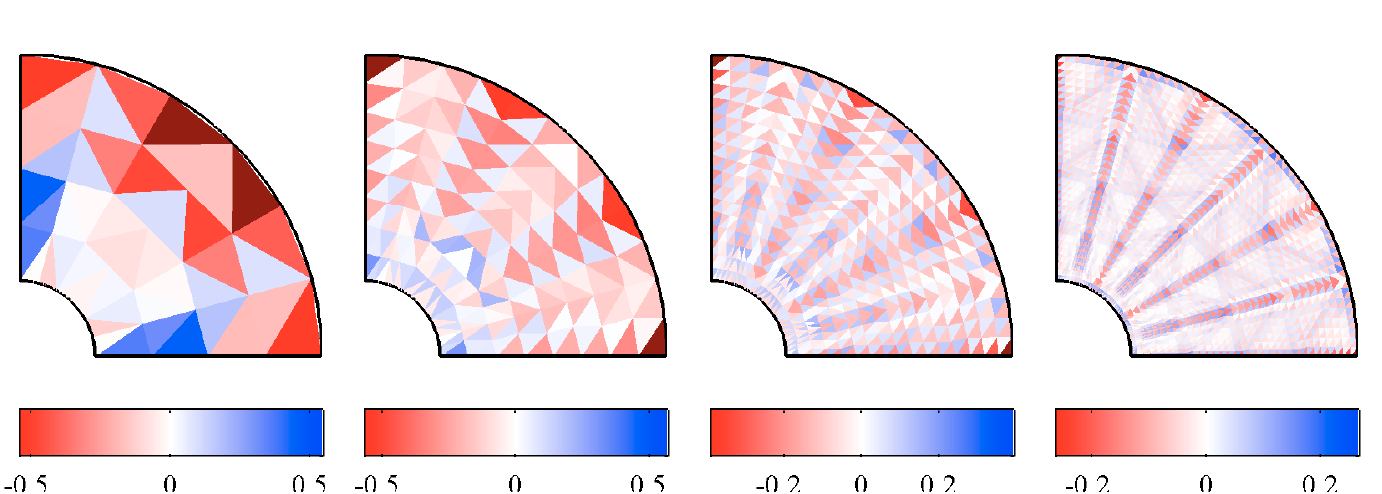}
	\caption{Thick-wall cylinder with TRI3 elements. Distribution of the effectivity index $D$ in uniformly refined meshes.}
	\label{fig:CylDTRI3}
\end{figure}

Figure~\ref{fig:CylDTRI3} shows the distribution of the local effectivity index $D$ in a sequence of meshes with linear triangular elements. The figure displays a quite uniform distribution of the local effectivity at each mesh, within the range between $[-0.54, 0.43 ]$ for the coarsest mesh. It is worth noting that the local values improve as we refine the meshes and that, for the last mesh in the sequence, the local effectivity is now within the range $[-0.26, 0.17]$. Some radial patterns in the distribution of $D$ can be seen as we increase the number of dof, which are attributed to local mesh configurations.  

\begin{figure}[!ht]
	\centering
	\includegraphics{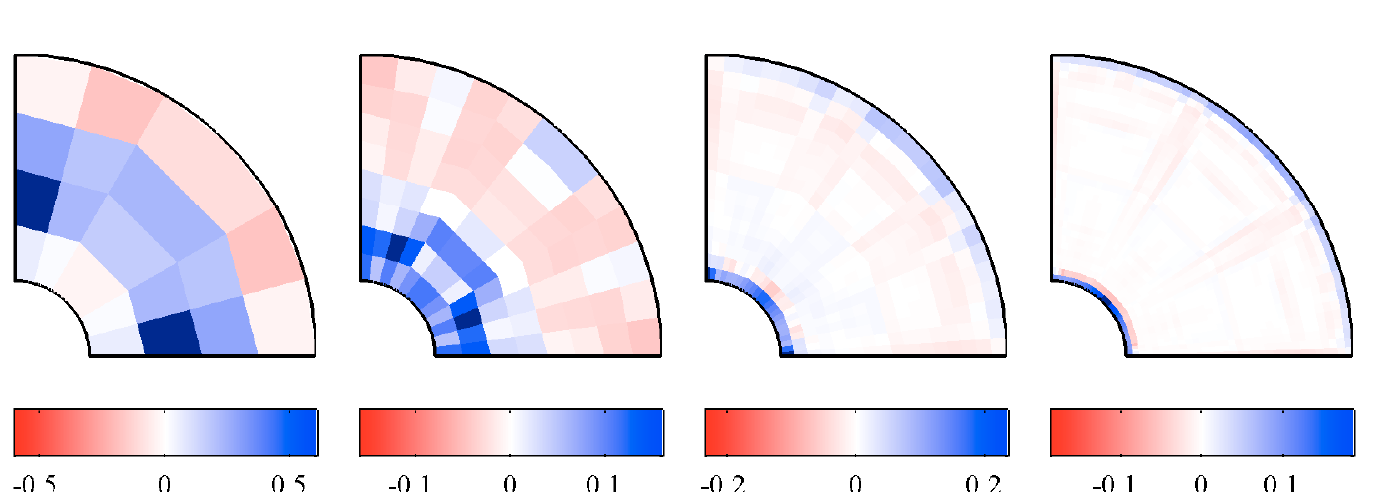}
	\caption{Thick-wall cylinder with QUAD4 elements. Distribution of the effectivity index $D$ in uniformly refined meshes.}
	\label{fig:CylDQUAD4}
\end{figure}

Figure~\ref{fig:CylDQUAD4} shows the distribution of the local effectivity index $D$ in a sequence of meshes with linear quadrilateral elements (QUAD4). Again, we can see a quite uniform distribution of the local effectivity at each mesh and that the local effectivity improves as the mesh is refined. This behaviour indicates that the proposed error estimator performs nicely at element level, which is important when guiding adaptive processes. 

\subsubsection{ Influence of support size.} 
 
One of the parameters that affects the performance of the proposed error estimator is the radius $R_{\Omega_{\vm{x}}}$ that defines the MLS support at a point $\vm{x}$, also known as the domain of influence. The idea is to define a support with a number of sampling points large enough to be able to solve the MLS fitting and to obtain an accurate polynomial expansion of the stresses, but not too large that we risk excessively smoothing the stress field and no longer describing the local behaviour of the solution. Moreover, larger supports means more computational effort as more sampling points should be considered.

In order to fix the domain of influence at a particular point we first evaluate the average size of the elements surrounding each node of the mesh. Then, we define the radius of the support at nodes as $R_{\Omega_{\vm{x}}}(\vm{x}_i) = k \, l(\vm{x}_i)$ where $k$ is a constant that takes positive values and $l(\vm{x}_i)$ is the average size of the elements containing node $i$. Once the value of $R_{\Omega_{\vm{x}}}$ is evaluated at nodes, the value of $R_{\Omega_{\vm{x}}}$ at any point $\vm{x}$ within an element is interpolated from the nodes using the displacement shape functions. Note that as $R_{\Omega_{\vm{x}}}$ is a function of $\vm{x}$. Its definition has been used in the derivatives of $s$ defined in (\ref{Eq:sFunction}) (or alternatively in (\ref{Eq:sFunctionMod})) required for the evaluation of (\ref{Eq:dWdx}, \ref{Eq:dWpdx}).

\begin{figure}[ht]
\centering
  \includegraphics{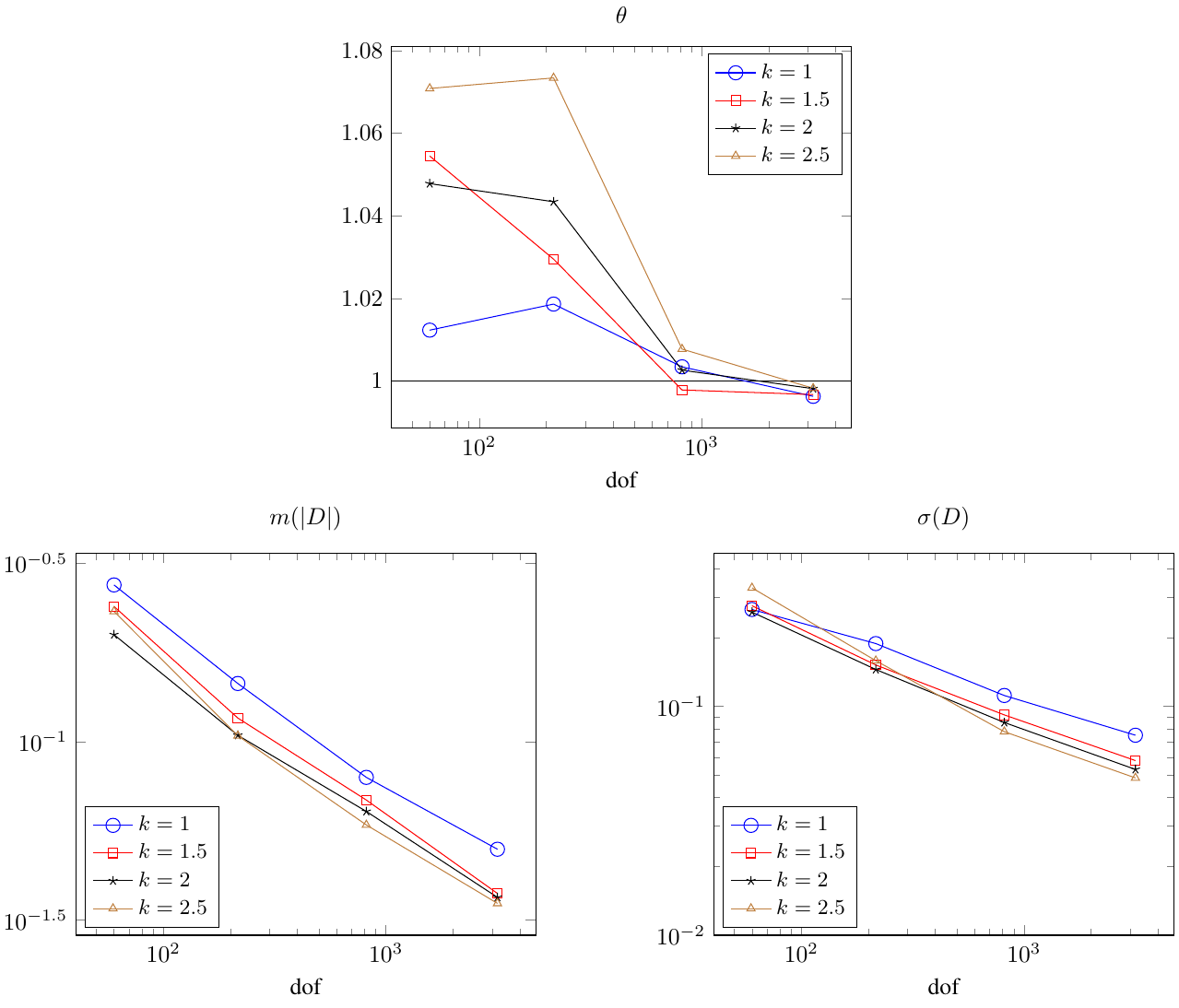}
  \caption{Thick-wall cylinder with uniformly refined TRI3 meshes. Global indicators $\theta$, $m(|D|)$ and $\sigma(D)$ for different values of $k$.}
  \label{fig:CylTRI3Ksize}
\end{figure}

\begin{figure}[ht]
\centering
    \includegraphics{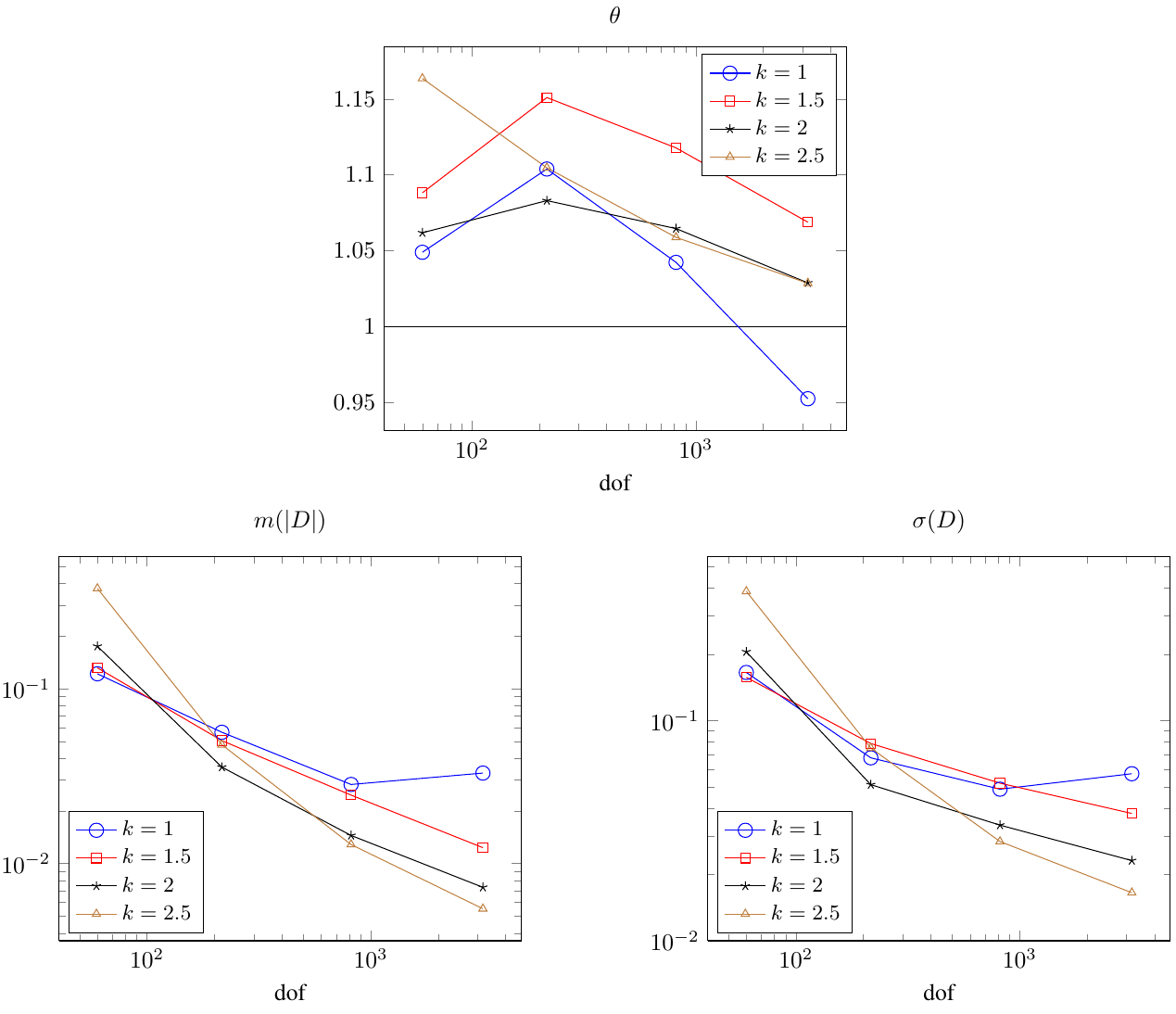}
  \caption{Thick-wall cylinder with uniformly refined QUAD4 meshes. Global indicators $\theta$, $m(|D|)$ and $\sigma(D)$ for different values of $k$.}
  \label{fig:CylQUAD4Ksize}
\end{figure}

Figure~\ref{fig:CylTRI3Ksize} shows the global results for a sequence of TRI3 elements considering different values of $k$. Figure~\ref{fig:CylQUAD4Ksize} shows the same results for QUAD4 meshes. Note that for small supports the values of the local indicators are less accurate even if the values for the global indicator are closer to one. A good balance between accuracy and local definition of the smoothing function is obtained for $k=2$, which is the value considered in the examples presented herein.

\subsection{ Westergaard problem -- FEM solution. } 
 
\begin{figure}[ht]
	\centering
	\includegraphics{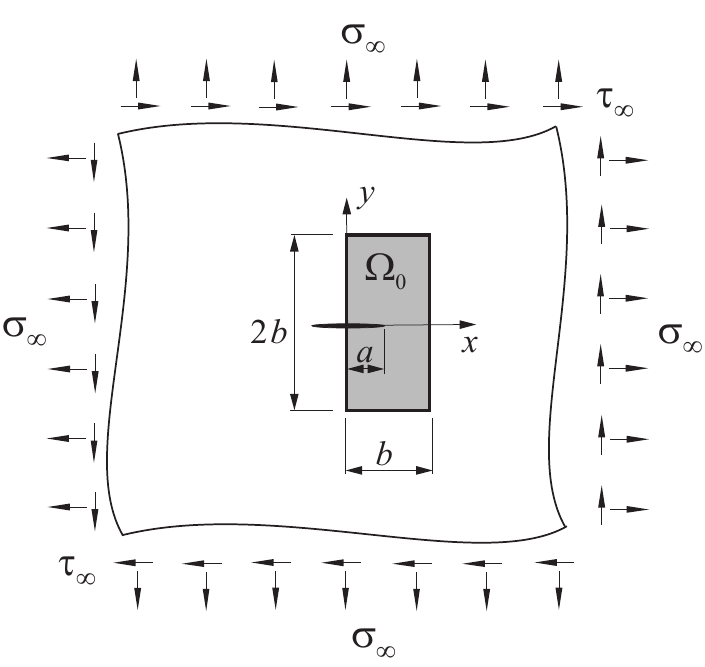}
	\caption{Westergaard problem. Infinite plate with a crack of length $2a$ under uniform tractions $\sigma_{\infty}$ (biaxial) and $\tau_{\infty}$. Finite portion of the domain $\Omega_0$, modelled with FE.}
	\label{fig:westergaard}
\end{figure}

To evaluate the performance of the proposed technique for singular problems we consider the Westergaard problem \cite{Rodenas2007a, Giner2005} as it has an exact analytical solution. The Westergaard problem corresponds to an infinite plate loaded at infinity with biaxial tractions $\sigma_{x \infty}=\sigma_{y \infty}=\sigma_{\infty}$ and shear traction $\tau_{\infty}$, presenting a crack of length $2a$ as shown in Figure~\ref{fig:westergaard}. Combining the externally applied loads we can obtain different loading conditions: pure mode I, pure mode II or mixed mode.  
 
The numerical model corresponds to a finite portion of the domain ($a=1$ and $b=4$ in Figure~\ref{fig:westergaard}). The applied projected stresses for mode I are evaluated from the analytical Westergaard solution \cite{Giner2005}:  

\begin{equation} \label{Eq:WesterStressI} 
\begin{array}{r@{\hspace{1ex}}c@{\hspace{1ex}}l} 
{\sigma _{x}^{I} } (x,y) & {=} & {\displaystyle \frac{\sigma _{\infty } }{\sqrt{\left|t\right|} } \bigg[\left(x\cos \frac{\phi }{2} -y\sin \frac{\phi }{2} \right)+y\frac{a^{2} }{\left|t\right|^{2} } \left(m\sin \frac{\phi }{2} -n\cos \frac{\phi }{2} \right)\bigg]} \\ 
\noalign{\medskip}{\sigma _{y}^{I} }(x,y) & {=} & {\displaystyle \frac{\sigma _{\infty } }{\sqrt{\left|t\right|} } \bigg[\left(x\cos \frac{\phi }{2} -y\sin \frac{\phi }{2} \right)-y\frac{a^{2} }{\left|t\right|^{2} } \left(m\sin \frac{\phi }{2} -n\cos \frac{\phi }{2} \right)\bigg]} \\ 
\noalign{\medskip}{ \tau _{xy}^{I} }(x,y) & {=} & {\displaystyle y\frac{a^{2} \sigma _{\infty } }{\left|t\right|^{2} \sqrt{\left|t\right|} } \left(m\cos \frac{\phi }{2} +n\sin \frac{\phi }{2} \right)} \end{array}
\end{equation}
\noindent and for mode II:
\begin{equation} \label{Eq:WesterStressII} 
\begin{array}{r@{\hspace{1ex}}c@{\hspace{1ex}}l} 
{\sigma _{x}^{II}}(x,y)  & {=} & {\displaystyle \frac{\tau _{\infty } }{\sqrt{\left|t\right|} } \bigg[2\left(y\cos \frac{\phi }{2} +x\sin \frac{\phi }{2} \right)-y\frac{a^{2} }{\left|t\right|^{2} } \left(m\cos \frac{\phi }{2} +n\sin \frac{\phi }{2} \right)\bigg]} \\ 
\noalign{\medskip}{\sigma _{y}^{II}}(x,y)  & {=} & {\displaystyle y\frac{a^{2} \tau _{\infty } }{\left|t\right|^{2} \sqrt{\left|t\right|} } \left(m\cos \frac{\phi }{2} +n\sin \frac{\phi }{2} \right)} \\ \noalign{\medskip}{\tau _{xy}^{II}}(x,y)  & {=} & {\displaystyle \frac{\tau _{\infty } }{\sqrt{\left|t\right|} } \bigg[\left(x\cos \frac{\phi }{2} -y\sin \frac{\phi }{2} \right)+y\frac{a^{2} }{\left|t\right|^{2} } \left(m\sin \frac{\phi }{2} -n\cos \frac{\phi }{2} \right)\bigg]} \end{array}
\end{equation}

\noindent where the stress fields are expressed as a function of $x$ and $y$, with origin at the centre of the crack. The parameters $t$, $m$, $n$ and $\phi$ are defined as
\begin{equation}
\begin{split}t& =(x+iy)^{2} -a^{2} =(x^{2} -y^{2} -a^{2} )+i(2xy)=m+in \\  m & =\textrm{Re}(t) =\textrm{Re}(z^{2} -a^{2} )=x^{2} -y^{2} -a^{2} \\ n & =\textrm{Im}(t)=(z^{2} -a^{2} )=2xy \\ \phi & =\textrm{Arg} (\bar{t})=\textrm{Arg} (m-in) \qquad\textrm{with }\phi \in \left[-\pi ,\pi \right], \; i^2=-1 \end{split} 
\end{equation}

For the problem analysed, the exact value of the SIF is given by 
\begin{equation} \label{Eq:SIFWestergaard} 
K_{I,ex} =\sigma_{\infty } \sqrt{\pi a} \qquad \qquad K_{II,ex} =\tau_{\infty } \sqrt{\pi a}  
\end{equation}

Material parameters are Young's modulus $E = 10^7$ and Poisson's ratio $ \nu= 0.333$. We consider loading conditions in pure mode I with $\sigma_{\infty} =100$ and $\tau_{\infty}=0$, pure  mode II with $\sigma_{\infty} =0$ and $\tau_{\infty}=100$, and mixed mode with $\sigma_{\infty} =100$ and $\tau_{\infty}=100$. To evaluate the SIF needed for recovering the singular part we use an equivalent domain integral technique, with a plateau function with radius $r_q = 0.9$ for the extraction \cite{Shih1986}. 

\begin{figure}[!ht]
  \centering
  \includegraphics{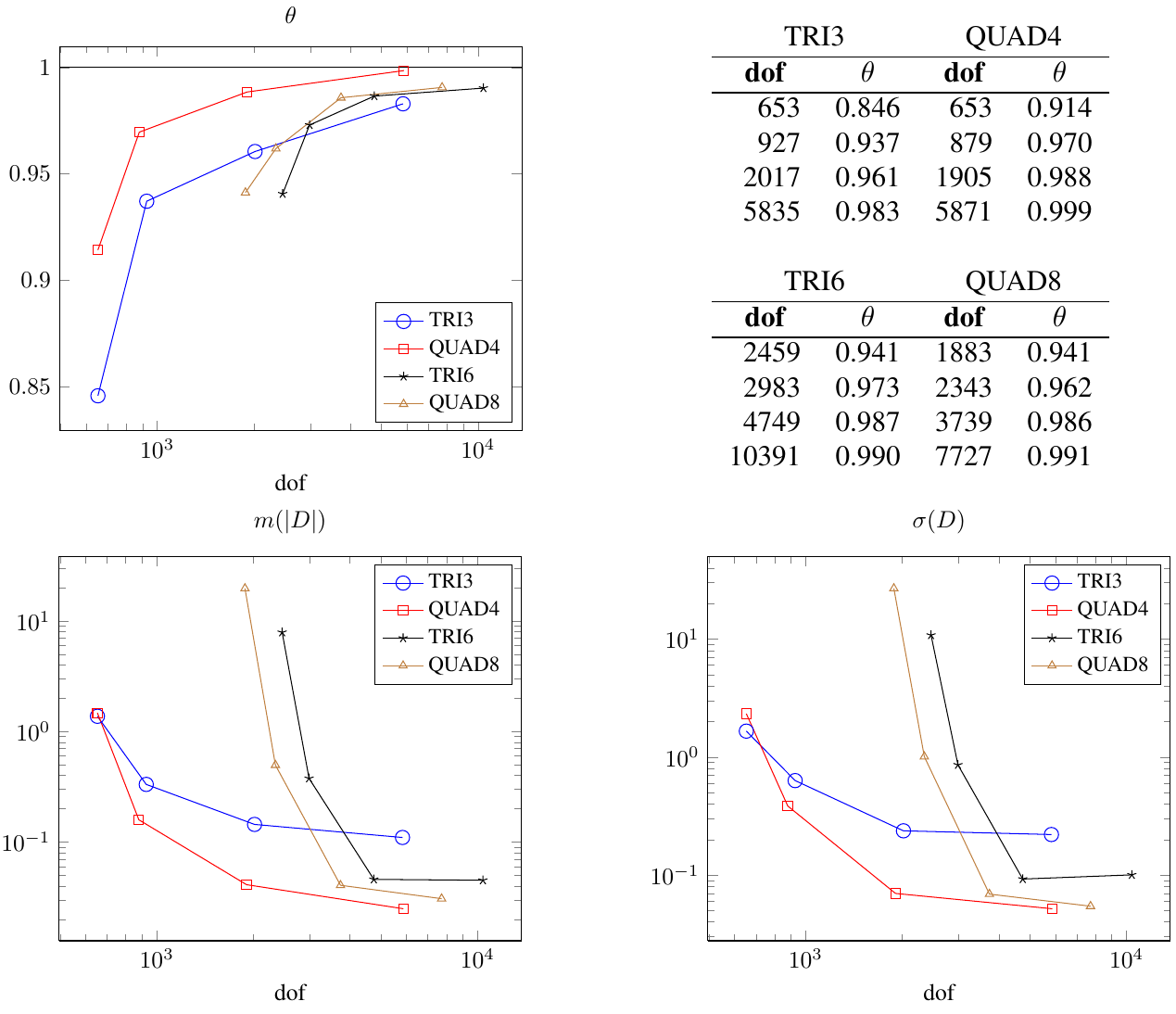}
  \caption{Westergaard problem under mode I with FEM \emph{h}-adapted meshes. Global indicators $\theta$, $m(|D|)$ and $\sigma(D)$.}
  \label{fig:FEMCrackEffectivity}
\end{figure}

Figure~\ref{fig:FEMCrackEffectivity} shows the evolution with respect to mesh refinement of the global parameters $\theta$, $m(|D|)$ and $\sigma(D)$  for different element types. In the figure, the global effectivity converges to the theoretical value of $\theta=1$ and both $m(|D|)$ and $\sigma(D)$ decrease with an increase of the number of dof. The performance of the proposed technique indicates an accurate error estimation for the meshes analysed. 

\begin{figure}[!ht]
	\centering
	\includegraphics{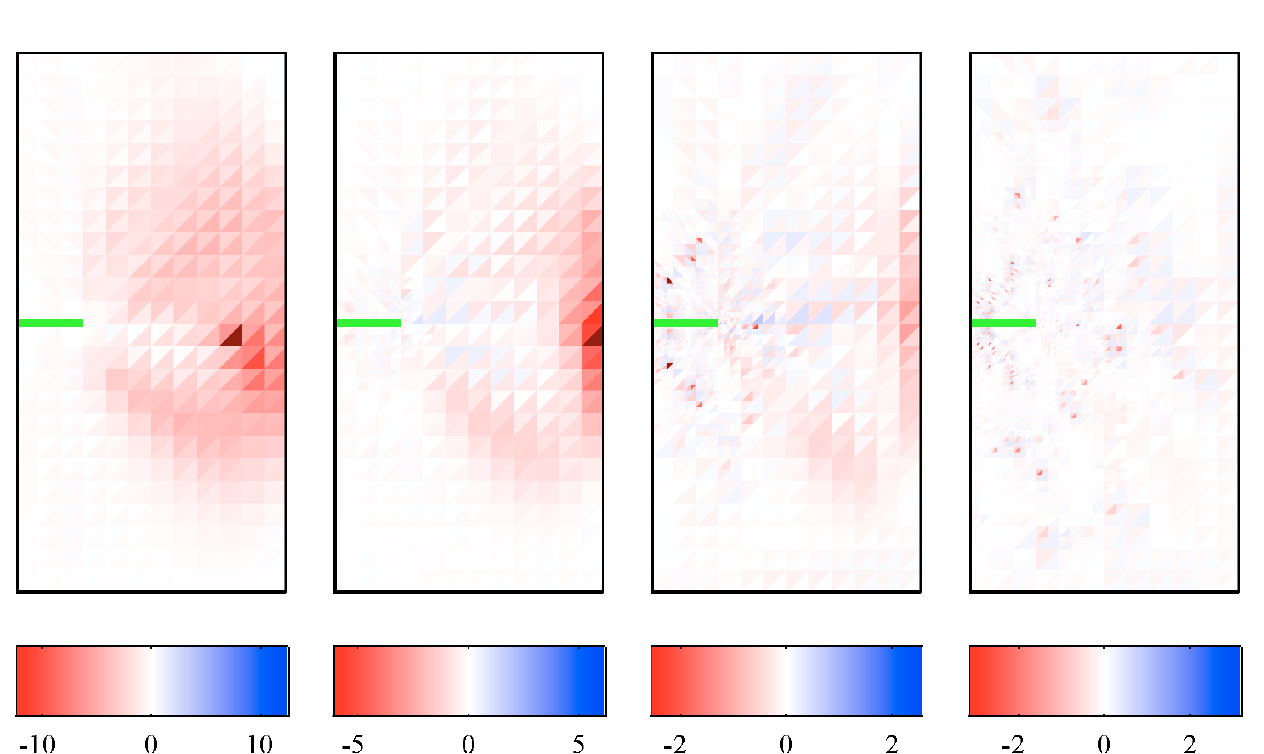}
	\caption{Westergaard problem under mode I with FEM and \emph{h}-adapted meshes of TRI3. Distribution of the effectivity index $D$.}
	\label{fig:FEMCRACKDTRI3}
\end{figure}

\begin{figure}[!ht]
	\centering
	\includegraphics{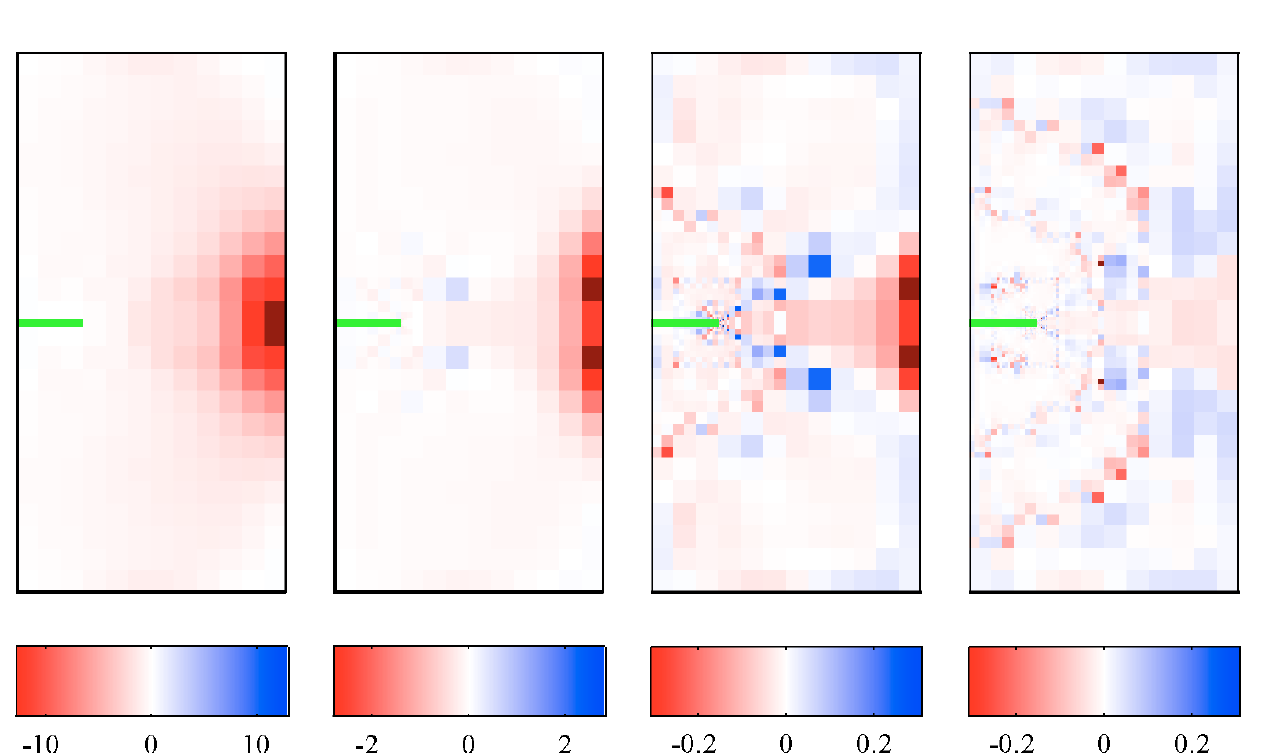}
	\caption{Westergaard problem under mode I with FEM and \emph{h}-adapted meshes of QUAD4. Distribution of the effectivity index $D$.}
	\label{fig:FEMCRACKDQUAD4}
\end{figure}

Figures~\ref{fig:FEMCRACKDTRI3} and \ref{fig:FEMCRACKDQUAD4} show the distribution of the local effectivity index $D$ in a sequence of TRI3 and QUAD4 meshes respectively. The splitting of the stress field into singular and smooth parts helps to recover a highly accurate stress field in the vicinity of the singular point. The distribution of the local effectivity index is homogeneous within the mesh and the values for $D$ decrease as we refine.

Because local error estimation techniques cannot take into account the pollution error due to the singularity, we can notice areas of the domain where the error is underestimated in this example, especially in the first meshes. The effect of pollution error is partially overcome by the use of \textit{h}-adaptive refinement (or enriched meshes as it is shown in the next section).

\subsection{ Westergaard problem -- XFEM solution.} 
 
Let us now consider the Westergaard problem from the previous section, solved using an enriched finite element approximation. In the numerical analyses, we use  a geometrical enrichment defined by a circular fixed enrichment area $B(x_0, r_e)$ with radius  $r_e = 0.5$, with its centre at the crack tip $x_0$ as proposed in \cite{Bechet2005}. For the extraction of the SIF we define a plateau function with radius $r_q = 0.9$ as in the FEM case. Bilinear elements are considered in the models. For the numerical integration of standard elements we use a $2\times2$ Gaussian quadrature rule. The elements intersected by the crack are split into triangular integration subdomains that do not contain the crack. Alternatives which do not require this subdivision are proposed in \cite{Ventura2006, Natarajan2010}. We use 7 Gauss points in each triangular subdomain, and a $5\times5$ quasipolar integration in the subdomains of the element containing the crack tip \cite{Bechet2005}, see Figure~\ref{fig:NodalSets}. We do not consider correction for blending elements. Methods to address blending errors are proposed in \cite{Chessa2003,Gracie2008,Fries2008,Tarancon2009}.

\begin{figure}[ht]
      \centering
  \includegraphics{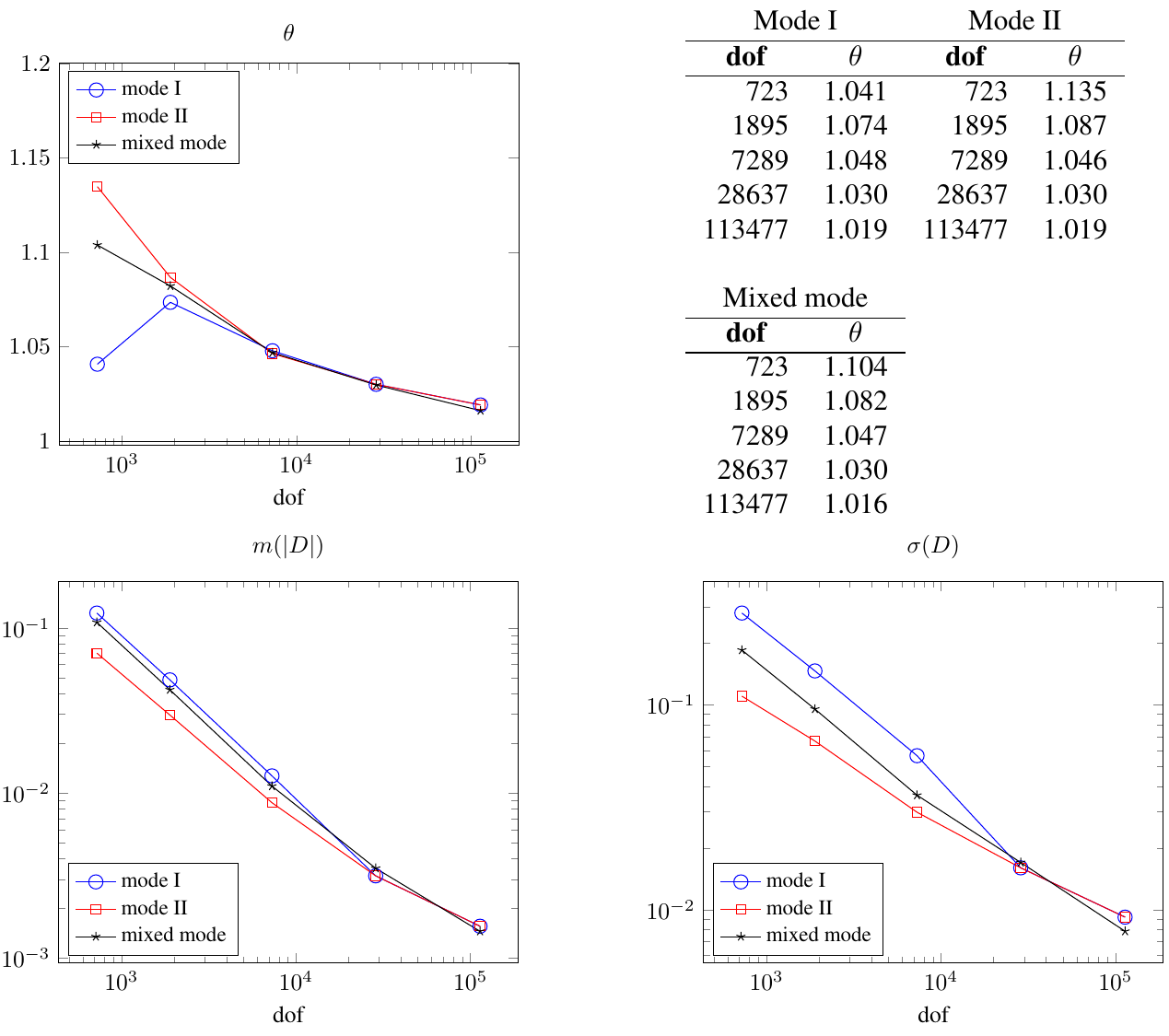}
  \caption{Westergaard problem with XFEM and structured meshes of QUAD4. Global indicators $\theta$, $m(|D|)$ and $\sigma(D)$.}
  \label{fig:CrackEffectivity}
\end{figure}

Figure~\ref{fig:CrackEffectivity} shows the evolution with respect to mesh refinement of the global parameters $\theta$, $m(|D|)$ and $\sigma(D)$  for the structured meshes of enriched QUAD4 elements. The curves represent the values obtained for the Westergaard problem under mode I, mode II and mixed mode loading conditions. In the figure, the global effectivity converges to the theoretical value of $\theta=1$ and both $m(|D|)$ and $\sigma(D)$ decrease with an increase of the number of dof. The results 
show that the proposed technique provides a sharp estimate of the true error.  

\begin{figure}[!ht]
	\centering
	\includegraphics{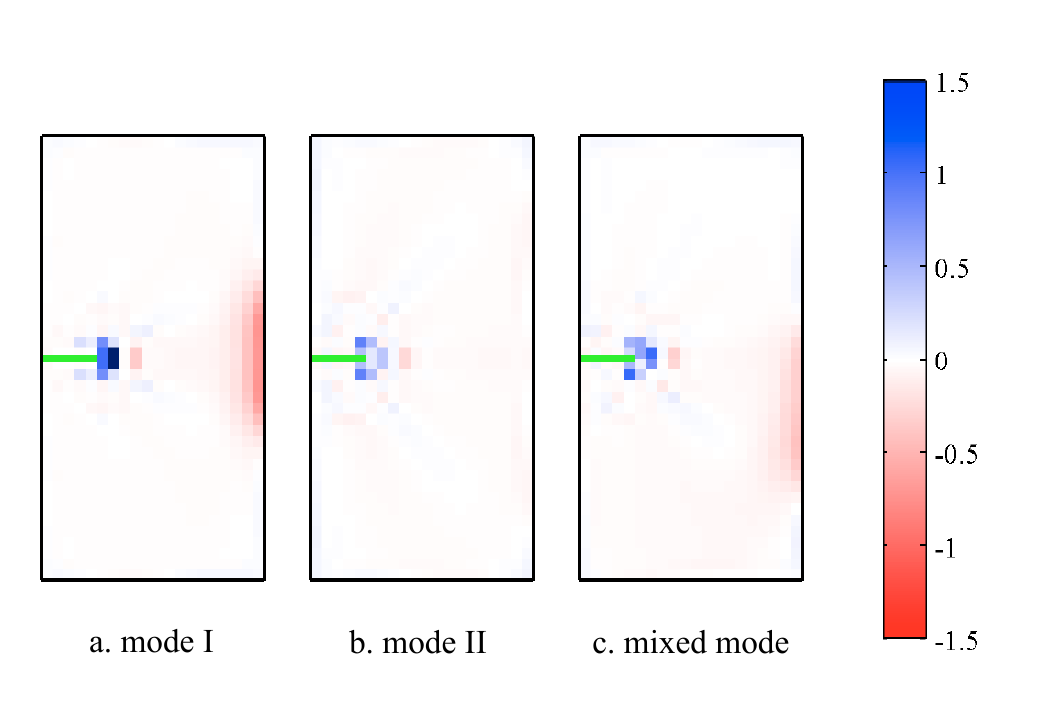}
	\caption{Westergaard problem with XFEM and structured meshes of QUAD4: a) mode I, b) mode II and c) mixed mode. Distribution of the effectivity index $D$.}
	\label{fig:CRACKDSM}
\end{figure}

Figure~\ref{fig:CRACKDSM} shows the distribution of $D$ in the second mesh (1895 dof) of the sequence of structured meshes for all the three loading modes. The results indicate a quite uniform distribution of the local effectivity. The values of $D$ indicate that the error at element level is accurately evaluated even where standard recovery techniques would produce the worst results: along the Neumann boundary, the crack faces and around the crack tip.

\begin{figure}[ht]
      \centering
  	\includegraphics{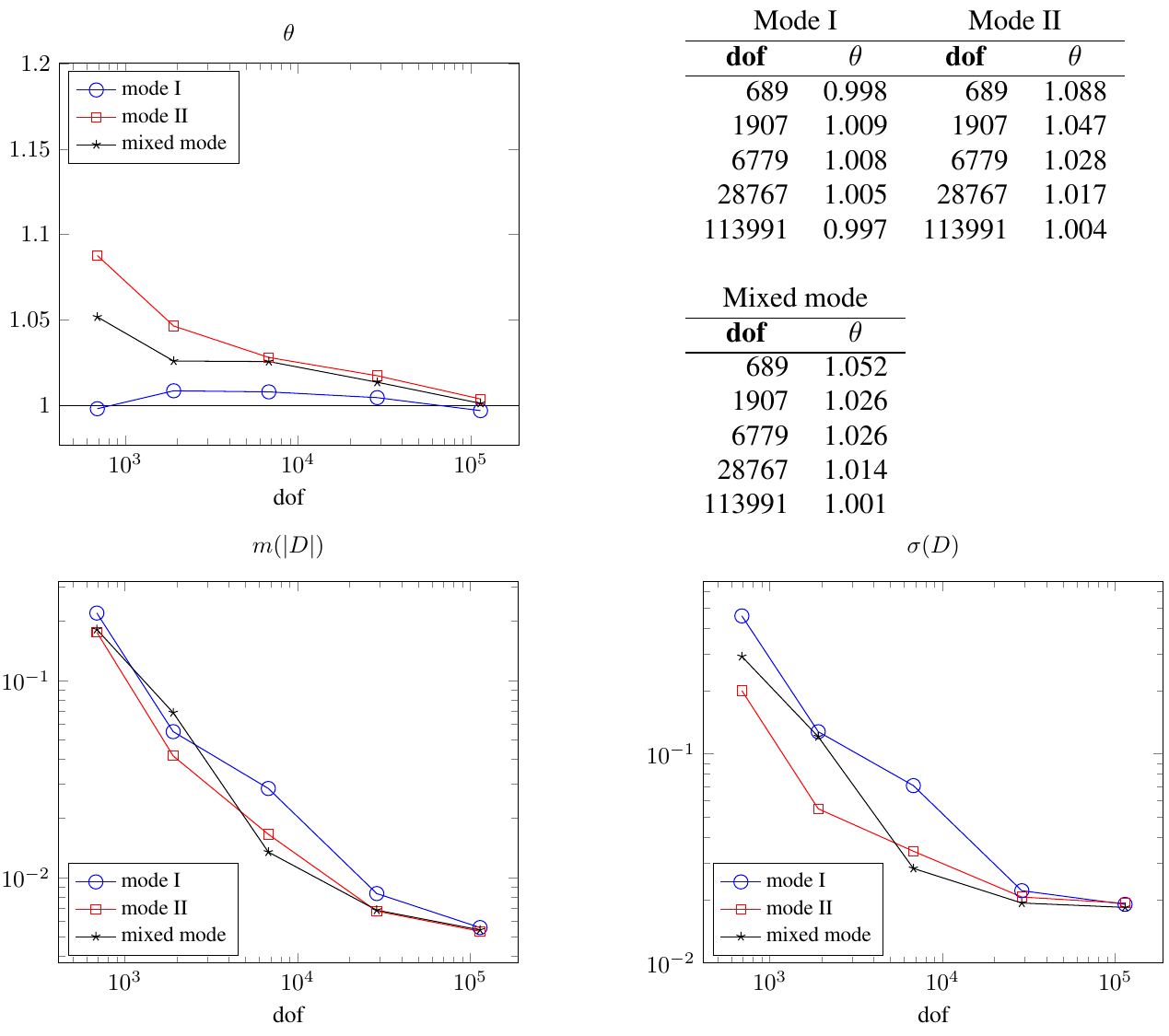}
  \caption{Westergaard problem with XFEM and non structured meshes of QUAD4. Global indicators $\theta$, $m(|D|)$ and $\sigma(D)$.}
  \label{fig:CrackEffectivityNS}
\end{figure}

\begin{figure}[!ht]
	\centering
	\includegraphics{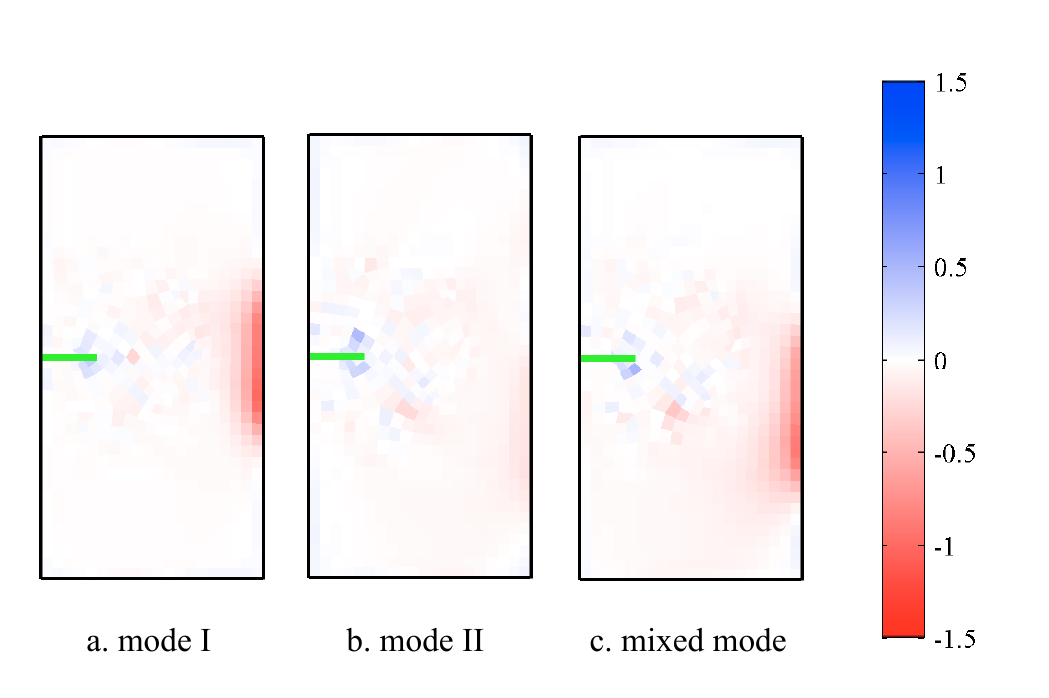}
	\caption{Westergaard problem with XFEM and non structured meshes of QUAD4: a) mode I, b) mode II and c) mixed mode. Distribution of the effectivity index $D$.}
	\label{fig:CRACKDNSM}
\end{figure}

For the case of non structured meshes the results for the same global parameters previously considered  are shown in Figure~\ref{fig:CrackEffectivityNS}. The local effectivity at element level for this meshes is depicted in Figure~\ref{fig:CRACKDNSM}. There is a similar behaviour to that seen for structured meshes. In general, the proposed technique exhibits an excellent performance when used to estimate the error in the XFEM approximations analysed.   

\begin{figure}[ht]
      \centering
  	\includegraphics{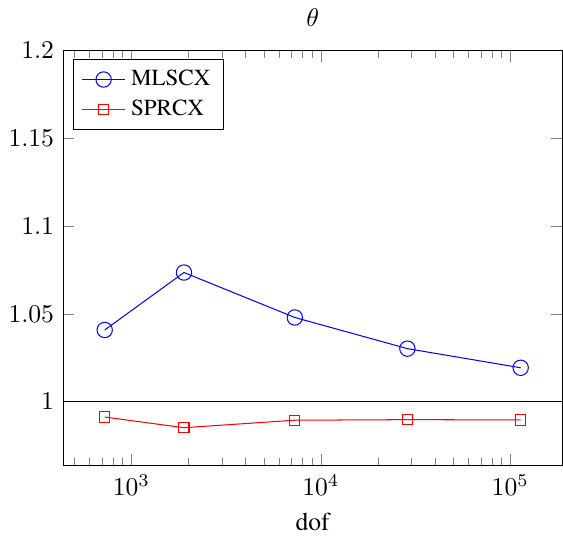}
  	\includegraphics{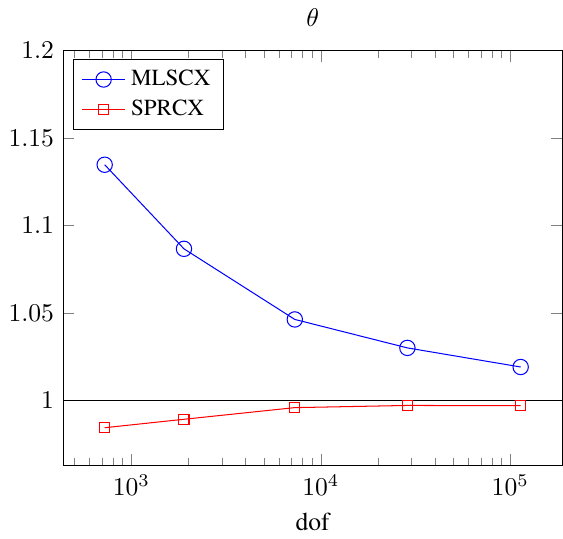}
  	\includegraphics{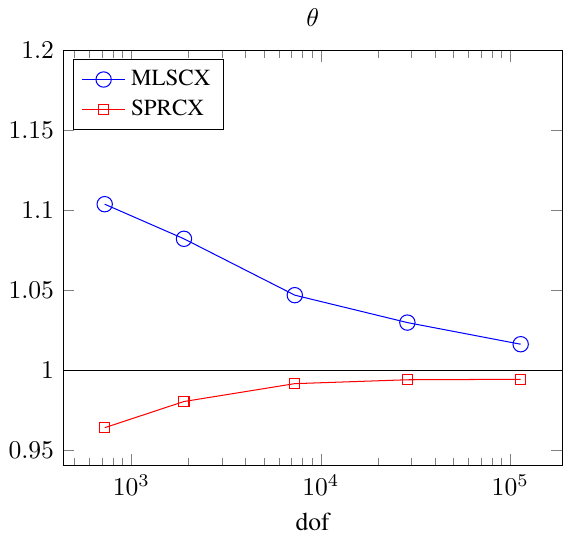}
  \caption{Westergaard problem with XFEM and structured meshes of QUAD4. Effectivity index $\theta$ for the MLSCX and SPRCX recovery techniques.}
  \label{fig:CrackSPRCX_Effectivity}
\end{figure}

\red{In Figure~\ref{fig:CrackSPRCX_Effectivity} we compare the results of the MLSCX with those of the SPRCX recovery procedure. In this case both techniques give values in the same order of magnitude}. 


\section{Conclusions}
\label{sec:Conclusions}

In this paper, the use of an equilibrated moving least squares recovery technique for FEM and XFEM problems has been investigated. The proposed technique uses a MLS approach to provide a continuous recovered stress field that enforces boundary equilibrium constraints. It also imposes a very accurate satisfaction, although not fully exact, of the internal equilibrium equation. Moreover, for singular problems it decomposes the stress field into two different parts, singular and smooth, in order to enable the technique to describe the singular behaviour of the solution. A visibility criterion is used near reentrant corners and cracks to properly define the weight of the sampling points within the support. 
 
The technique presented here has been validated using four different examples with known analytical solution. The numerical results have shown the accuracy of the proposed technique, which provides values of the effectivity index that converge and are very close to the theoretical value $\theta=1$. The distribution of the local effectivity at the elements is homogeneous for the tests considered, and the mean value $m(|D|)$ and standard deviation $\sigma(D)$ decrease as we increase the number of dof. \red{The obtained MLS recovered field is not fully statically admissible, thus, the procedure does not guarantee the upper bound property. For this reason, it nearly bounds the exact error but not always yields an effectivity index greater than one, as clearly seen in the first example.} In any case, the numerical results show that for the examples presented the proposed technique yields sharp error estimates, which are very accurate when compared with previous MLS approaches. \red{Extension of this work to 3D problems is feasible given that the SIF along the crack front is evaluated with sufficient accuracy. It is known that in 3D problems the evaluation of the SIF is less accurate. In \cite{Rodenas2007a} the influence of the accuracy in the evaluation of the SIF in the error estimator is investigated.}

 
\section{Acknowledgements}

This work has been carried within the framework of the research project DPI2010-20542 of the Ministerio de Ciencia y e Innovación y  (Spain). The financial support of the Universitat Politècnica de València and Generalitat Valenciana is also acknowledged.\\
Support from the EPSRC grant EP/G042705/1 ``Increased Reliability for Industrially Relevant Automatic Crack Growth Simulation with the eXtended Finite Element Method'' is acknowledged.
  




\bibliographystyle{bibstyle}
\bibliography{General}

\end{document}